\documentclass[leqno,12pt]{amsart}
\usepackage{amssymb} 
\theoremstyle{plain} % ??E???
\newtheorem{theorem}{\indent\sc Theorem}[section] % ???I?A?E?A??

\newtheorem{proposition}[theorem]{\indent\sc Proposition}

\theoremstyle{definition} % I??O?E?a

\begin{document}

\title[Polynomial Hamiltonian system in two variables]{Polynomial Hamiltonian system in two variables with $W({A}^{(1)}_1)$-symmetry and the second Painlev\'e hierarchy \\}

\author{By\\
Yusuke Sasano}

\renewcommand{\thefootnote}{\fnsymbol{footnote}}
\footnote[0]{2000\textit{ Mathematics Subjet Classification}.
34M55; 34M45; 58F05; 32S65.}

%\date{August, 1, 2007}

\begin{abstract} We find a one-parameter family of polynomial Hamiltonian system in two variables with $W({A}^{(1)}_1)$-symmetry. We also show that this system can be obtained by the compatibility conditions for the linear differential equations in three variables. We give a relation between it and the second member of the second Painlev\'e hierarchy. Moreover, we give some relations between an autonomous version of its polynomial Hamiltonian system in two variables and the mKdV hierarchies.
\end{abstract}

\maketitle

\section{Introduction}
In this paper, we find a 1-parameter family of total differential system in two variables $t,s$
\begin{equation}\label{eq:1}
  \left\{
  \begin{aligned}
   dx =&(-2xp-\alpha_0)dt+(-2xz-2\alpha_0 w)ds,\\
   dy =&(2yp+\alpha_1)dt+(2yz-2\alpha_1 q)ds,\\
   dz =&\left(\frac{x+y}{2} \right)dt+(-xq+yw-2p)ds,\\
   dw =&(z-2wp)dt+(-2xp-2\alpha_0-\alpha_1)ds,\\
   dq =&(z+2qp)dt+(2yp+\alpha_0+2\alpha_1)ds,\\
   dp =&\left(\frac{w+q}{2} \right)dt+\left(\frac{x+y}{2} \right)ds,
   \end{aligned}
  \right.
\end{equation}
where the constant parameters $\alpha_i$ satisfy the relation $\alpha_0+\alpha_1=1$.

In next theorem, we show that this system satisfies the compatibility conditions, and we can characterize this system by the following holomorphy conditions $r_0,r_1$.
\begin{theorem}\label{th:1}
Let us consider a system of first order ordinary differential equations in the polynomial class\rm{:\rm}
\begin{equation*}
\frac{dx}{dT}=f_1(x,y,z,w,q,p), \ldots, \frac{dp}{dT}=f_6(x,y,z,w,q,p) \quad f_i \in {\mathbb C}[x,y,z,w,q,p] \quad (i=1,\ldots,6).
\end{equation*}
We assume that

$(A1)$ $deg(f_i)=2$ with respect to $x,y,z,w,q,p$.

$(A2)$ The right-hand side of this system becomes again a polynomial in each coordinate system $r_i \ (i=0,1):$
\begin{align*}
\begin{split}
&r_0:x_0=-(xp+\alpha_0)p, \quad y_0=y-4zp+4wp^2, \quad z_0=z-2wp, \quad w_0=w,\\
&q_0=q-2p^2, \quad p_0=\frac{1}{p},\\
&r_1:x_1=x+4zp+4qp^2, \quad y_1=-(yp+\alpha_1)p, \quad z_1=z+2qp, \quad w_1=w+2p^2,\\
&q_1=q, \quad p_1=\frac{1}{p}.
\end{split}
\end{align*}
Then we can obtain two systems:
\begin{equation}
  \left\{
  \begin{aligned}
   \frac{\partial x}{\partial t} =&-2xp-\alpha_0,\\
   \frac{\partial y}{\partial t} =&2yp+\alpha_1,\\
   \frac{\partial z}{\partial t} =&\frac{x+y}{2},\\
   \frac{\partial w}{\partial t} =&z-2wp,\\
   \frac{\partial q}{\partial t} =&z+2qp,\\
   \frac{\partial p}{\partial t} =&\frac{w+q}{2}
   \end{aligned}
  \right.
\end{equation}
and
\begin{equation}
  \left\{
  \begin{aligned}
   \frac{\partial x}{\partial s} =&-2xz-2\alpha_0 w,\\
   \frac{\partial y}{\partial s} =&2yz-2\alpha_1 q,\\
   \frac{\partial z}{\partial s} =&-xq+yw-2p,\\
   \frac{\partial w}{\partial s} =&-2xp-2\alpha_0-\alpha_1,\\
   \frac{\partial q}{\partial s} =&2yp+\alpha_0+2\alpha_1,\\
   \frac{\partial p}{\partial s} =&\frac{x+y}{2}.
   \end{aligned}
  \right.
\end{equation}
These two systems satisfy the compatibility conditions:
\begin{equation}
\frac{\partial }{\partial s} \frac{\partial x}{\partial t}=\frac{\partial }{\partial t} \frac{\partial x}{\partial s}, \quad \frac{\partial }{\partial s} \frac{\partial y}{\partial t}=\frac{\partial }{\partial t} \frac{\partial y}{\partial s}, \ldots, \frac{\partial }{\partial s} \frac{\partial p}{\partial t}=\frac{\partial }{\partial t} \frac{\partial p}{\partial s}.
\end{equation}
\end{theorem}
We remark that these transition functions in $r_0,r_1$ satisfy the condition:
\begin{equation}
dx_i \wedge dy_i \wedge dz_i \wedge dw_i \wedge dq_i \wedge dp_i=dx \wedge dy \wedge dz \wedge dw \wedge dq \wedge dp \quad (i=0,1).
\end{equation}

\begin{theorem}\label{th:2}
The system \eqref{eq:1} admits the extended affine Weyl group symmetry of type $A_1^{(1)}$ as the group of its B{\"a}cklund transformations, whose generators $s_0,s_1,{\pi}$ defined as follows$:$ with {\it the notation} $(*):=(x,y,z,w,q,p;\alpha_0,\alpha_1)$\rm{: \rm}
\begin{align*}
s_0:(*) \rightarrow &\left(x,y+\frac{4\alpha_0 z}{x}+\frac{4\alpha_0^2 w}{x^2},z+\frac{2\alpha_0 w}{x},w,q+\frac{4\alpha_0 p}{x}+\frac{2\alpha_0^2}{x^2},p+\frac{\alpha_0}{x};-\alpha_0,\alpha_1+2\alpha_0 \right),\\
s_1:(*) \rightarrow &\left(x-\frac{4\alpha_1 z}{y}+\frac{4\alpha_1^2 q}{y^2},y,z-\frac{2\alpha_1 q}{y},w-\frac{4\alpha_1 p}{y}-\frac{2\alpha_1^2}{y^2},q,p+\frac{\alpha_1}{y};\alpha_0+2\alpha_1,-\alpha_1 \right),\\
\pi:(*) \rightarrow &(-y,-x,-z,-q,-w,-p;\alpha_1,\alpha_0).
\end{align*}
\end{theorem}

\section{Lax pair for the system}
In this section, we show that the system \eqref{eq:1} can be obtained by the compatibility conditions for the linear differential equations of second order in three variables $T,t,s$
\begin{equation}
T{\partial}_{T} \vec{u}=A(T,t)\vec{u}, \quad {\partial}_{t} \vec{u}=B_1(T,t)\vec{u}, \quad \vec{u}=\begin{pmatrix}
u_1\\
u_2 
\end{pmatrix},
\end{equation}
and
\begin{equation}
T{\partial}_{T} \vec{u}=A(T,s)\vec{u}, \quad {\partial}_{s} \vec{u}=B_3(T,s)\vec{u}.
\end{equation}

The matrices $A(T,t),B_1(T,t)$ and $B_3(T,s)$ are explicitly given by
\begin{align}
\begin{split}
&-A(T,t)=\begin{pmatrix}
\varepsilon_1 & y\\
0 & \varepsilon_2 
\end{pmatrix}+\begin{pmatrix}
2z & 4q\\
-x & -2z 
\end{pmatrix}T+\begin{pmatrix}
8p & -8\\
-4w & -8p
\end{pmatrix}T^2+\begin{pmatrix}
0 & 0\\
-8 & 0 
\end{pmatrix}T^3,\\
&B_1(T,t)=\begin{pmatrix}
p & -1\\
0 & -p 
\end{pmatrix}+\begin{pmatrix}
0 & 0\\
-1 & 0
\end{pmatrix}T,\\
&B_3(T,s)=\begin{pmatrix}
\frac{z}{4} & \frac{q}{2}\\
0 & -\frac{z}{4} 
\end{pmatrix}+\begin{pmatrix}
p & -1\\
-\frac{w}{2} & -p 
\end{pmatrix}T+\begin{pmatrix}
0 & 0\\
-1 & 0
\end{pmatrix}T^2.
\end{split}
\end{align}
Here, the matrices $A(T,t),B_1(T,t)$ and $B_3(T,s)$ depend on $t,s$, and $\varepsilon_i$ are constant parameters.

At first, by the compatibility conditions
\begin{equation}
{\partial}_{t}(A(T,t))-T{\partial}_{T}(B_1(T,t))+[A(T,t),B_1(T,t)]=0
\end{equation}
for the linear differential equations in two variables $T,t$
\begin{align}
\begin{split}
&T{\partial}_{T} \vec{u}=A(T,t)\vec{u}, \quad {\partial}_{t} \vec{u}=B_1(T,t)\vec{u},
\end{split}
\end{align}
we can obtain the non-linear ordinary differential system in the variable $t$
\begin{equation}
  \left\{
  \begin{aligned}
   \frac{dx}{dt} &=-2xp-1+\varepsilon_1-\varepsilon_2,\\
   \frac{dy}{dt} &=2yp+\varepsilon_1-\varepsilon_2,\\
   \frac{dz}{dt} &=\frac{x+y}{2},\\
   \frac{dw}{dt} &=-2wp+z,\\
   \frac{dq}{dt} &=2qp+z,\\
   \frac{dp}{dt} &=\frac{w+q}{2},
   \end{aligned}
  \right. 
\end{equation}
where $\varepsilon_i$ satisfy the relations
\begin{equation}
  \left\{
  \begin{aligned}
  \alpha_0  &=-\varepsilon_1+\varepsilon_2+1,\\
  \alpha_1  &=\varepsilon_1-\varepsilon_2.
   \end{aligned}
  \right. 
\end{equation}
Next, by the compatibility conditions
\begin{equation}
{\partial}_{s}(A(T,s))-T{\partial}_{T}(B_3(T,s))+[A(T,s),B_3(T,s)]=0
\end{equation}
for the linear differential equations in two variables $T,s$
\begin{align}
\begin{split}
&T{\partial}_{T} \vec{u}=A(T,s)\vec{u}, \quad {\partial}_{s} \vec{u}=B_3(T,s)\vec{u},
\end{split}
\end{align}
we can obtain the non-linear ordinary differential system in the variable $s$
\begin{equation}
  \left\{
  \begin{aligned}
   \frac{dx}{ds} &=-2xz-2(1-\varepsilon_1+\varepsilon_2)w,\\
   \frac{dy}{ds} &=2yz-2(\varepsilon_1-\varepsilon_2)q,\\
   \frac{dz}{ds} &=-xq+yw-2p,\\
   \frac{dw}{ds} &=-2xp-2(1-\varepsilon_1+\varepsilon_2)-(\varepsilon_1-\varepsilon_2),\\
   \frac{dq}{ds} &=2yp+(1-\varepsilon_1+\varepsilon_2)+2(\varepsilon_1-\varepsilon_2),\\
   \frac{dp}{ds} &=\frac{x+y}{2}.
   \end{aligned}
  \right. 
\end{equation}

\section{Polynomial Hamiltonian system}
In this section, we show that the system \eqref{eq:1} is equivalent to the polynomial Hamiltonian system in two variables $t,s$. At first, we show that the system \eqref{eq:1} has two first integrals.
\begin{proposition}
This system \eqref{eq:1} has its first integrals:
\begin{equation}
  \left\{
  \begin{aligned}
    w-q+2p^2+3s=&C_1,\\
    8sp^2+4zp-2wq+x-y+4sw-4sq+t+6s^2=&C_2, \quad (C_1,C_2 \in {\mathbb C}).
   \end{aligned}
  \right.
\end{equation}
\end{proposition}

\begin{theorem}
The transformations
\begin{equation}\label{eq:2}
  \left\{
  \begin{aligned}
   q_1 =&p,\\
   p_1 =&x,\\
   q_2 =&z-2wp,\\
   p_2 =&w
   \end{aligned}
  \right. 
\end{equation}
take the system \eqref{eq:1} to the polynomial Hamiltonian system in two variables $t,s$
\begin{equation}\label{eq:3}
  \left\{
  \begin{aligned}
   dq_1 =&\frac{\partial H_1}{\partial p_1}dt+\frac{\partial H_2}{\partial p_1}ds,\\
   dp_1 =&-\frac{\partial H_1}{\partial q_1}dt-\frac{\partial H_2}{\partial q_1}ds,\\
   dq_2 =&\frac{\partial H_1}{\partial p_2}dt+\frac{\partial H_2}{\partial p_2}ds,\\
   dp_2 =&-\frac{\partial H_1}{\partial q_2}dt-\frac{\partial H_2}{\partial q_2}ds
   \end{aligned}
  \right. 
\end{equation}
with the polynomial Hamiltonians \rm{(cf. \cite{Mazzocco}) \rm}
\begin{align}
\begin{split}
H_1=&q_1^2p_1+\frac{3}{2}sp_1-\frac{C_1}{2}p_1+\alpha_0 q_1\\
&-p_2^3+(C_1-3s)p_2^2+\left(-3s^2+2C_1 s+\frac{t}{2}-\frac{C_2}{2} \right)p_2-\frac{q_2^2}{2}+p_1 p_2,\\
H_2=&\frac{p_1^2}{2}-\frac{1}{2}(C_2-4C_1 s+6s^2-t)p_1+(2\alpha_0+\alpha_1)q_2\\
&-p_1p_2^2+2q_1^2 p_1p_2+2q_1p_1q_2+2\alpha_0 q_1p_2+(C_1-3s)p_1p_2.
\end{split}
\end{align}
\end{theorem}
The relations between $x,y,z,w,q,p$ and $q_1,p_1,q_2,p_2$ are given by
\begin{equation}
  \left\{
  \begin{aligned}
   x =&p_1,\\
   y =&4q_1^2 p_2-2p_2^2+4q_1q_2+2C_1 p_2-6s p_2+p_1-6s^2+t+4C_1 s-C_2,\\
   z =&q_2+2q_1 p_2,\\
   w =&p_2,\\
   q =&2q_1^2+p_2+3s-C_1,\\
   p =&q_1.
   \end{aligned}
  \right.
\end{equation}

Setting $C_1=0,C_2=0,s=0$ and $\alpha_0:=\frac{1}{2}-\alpha_2$ in the Hamiltonian $H_1$, we can obtain the Hamiltonian in the variable $t$ given by \cite{Sasano1}. This Hamiltonian system is equivalent to the second member of the second Painlev\'e hierarchy (see \cite{Sasano1})
\begin{equation}
P_{II}^{(2)}:\frac{d^4u}{dt^4}=10u \left(\frac{du}{dt} \right)^2+10u^2 \frac{d^2u}{dt^2}-6u^5+tu+\alpha_2 \quad (\alpha_2 \in {\mathbb C}).
\end{equation}

\section{Symmetry and holomorphy  of the system \eqref{eq:3}}
\begin{theorem}
The system \eqref{eq:3} admits the extended affine Weyl group symmetry of type $A_1^{(1)}$ as the group of its B{\"a}cklund transformations, whose generators $s_0,s_1,{\pi}$ defined as follows$:$ with {\it the notation} $(*):=(q_1,p_1,q_2,p_2,t,s;\alpha_0,\alpha_0)$\rm{: \rm}
\begin{align*}
s_0:(*) \rightarrow &\left(q_1+\frac{\alpha_0}{f_0},p_1,q_2,p_2,t,s;-\alpha_0,\alpha_1+2\alpha_0 \right),\\
s_1:(*) \rightarrow &(q_1+\frac{\alpha_1}{f_1},p_1-\frac{4\alpha_1(q_2+2q_1p_2)}{f_1}+\frac{4\alpha_1^2(p_2+2q_1^2-(C_1-3s))}{f_1^2},\\
&q_2-\frac{2\alpha_1(2p_2-2q_1^2-(C_1-3s))}{f_1}+\frac{12\alpha_1^2 q_1}{f_1^2}+\frac{4\alpha_1^3}{f_1^3},p_2-\frac{4\alpha_1 q_1}{f_1}-\frac{2\alpha_1^2}{f_1^2},t,s;\\
&\alpha_0+2\alpha_1,-\alpha_1),\\
\pi:(*) \rightarrow &(-q_1,-f_1,-(q_2+4q_1(q_1^2+p_2)-2(C_1-3s)q_1),-(p_2+2q_1^2-C_1+3s),t,s;\alpha_1,\alpha_0),
\end{align*}
where $f_0:=p_1$ and $f_1:=p_1+4q_1^2p_2-2p_2^2+4q_1q_2+2(C_1-3s)p_2-C_2+4C_1s-6s^2+t$.
\end{theorem}
We note that the B{\"a}cklund transformations of this system satisfy the universal description for $A_1^{(1)}$ root system:
\begin{equation}
s_i(g)=g+\frac{\alpha_i}{f_i}\{f_i,g\}+\frac{1}{2!}\left(\frac{\alpha_i}{f_i} \right)^2 \{f_i,\{f_i,g\}\}+\cdots \quad (g \in {\mathbb C}(t,s)[q_1,p_1,q_2,p_2]),
\end{equation}
where $\{p_i,q_j\}=\delta_{ij}$ and $\{p_i,p_j\}=\{q_i,q_j\}=0$.

\begin{theorem}
Let us consider a polynomial Hamiltonian system with Hamiltonian $H \in {\mathbb C}(t,s)[q_1,p_1,q_2,p_2]$. We assume that

$(B1)$ $deg(H)=5$ with respect to $q_1,p_1,q_2,p_2$.

$(B2)$ This system becomes again a polynomial Hamiltonian system in each coordinate $\tilde{r}_i \ (i=0,1)${\rm : \rm}
\begin{align*}
\begin{split}
\tilde{r}_0:(x_0,y_0,z_0,w_0)=&\left(\frac{1}{q_1},-(q_1f_0+\alpha_0)q_1,q_2,p_2 \right),\\
\tilde{r}_1:(x_1,y_1,z_1,w_1)=&\left(\frac{1}{q_1},-(q_1f_1+\alpha_1)q_1,q_2+4q_1(q_1^2+p_2)-2(C_1-3s)q_1,p_2+2q_1^2 \right).
\end{split}
\end{align*}
Then such a system coincides with the Hamiltonian system \eqref{eq:3} with the polynomial Hamiltonians $H_1,H_2$.
\end{theorem}
We note that the conditions $(B2)$ should be read that
\begin{align*}
\begin{split}
&r_0(H), \quad r_0(H-q_1),\\
&r_1(H), \quad r_1(H+8q_1^3+6q_1p_2-4(C_1-3s)q_1)
\end{split}
\end{align*}
are polynomials with respect to $x_i,y_i,z_i,w_i$.

Next, we study a solution of the system \eqref{eq:3} which is written by the use of known functions.

By the transformation $\pi$, the fixed solution is derived from
\begin{align}
\begin{split}
&\alpha_0=\alpha_1, \quad q_1=-q_1, \quad p_1=-f_1,\\
&q_2=-(q_2+4q_1(q_1^2+p_2)-2(C_1-3s)q_1), \quad p_2=-(p_2+2q_1^2-C_1+3s).
\end{split}
\end{align}
Then we obtain
\begin{equation}
(q_1,p_1,q_2,p_2;\alpha_0,\alpha_1)=\left(0,\frac{-2t+3s^2-2C_1s-C_1^2+2C_2}{4},0,\frac{C_1-3s}{2};\frac{1}{2},\frac{1}{2} \right).
\end{equation}
Applying the B{\"a}cklund transformations to this seed solution, we can obtain a series of special solutions.

Finally, let us consider the relation between the polynomial Hamiltonian system \eqref{eq:3} and a modified type of the mKdV equation. In this paper, we can make the birational transformations between the polynomial Hamiltonian system \eqref{eq:3} and a modified type of mKdV equation.
\begin{theorem}
The birational transformations
\begin{equation}\label{eq:111111}
  \left\{
  \begin{aligned}
   x =&q_1,\\
   y =&q_1^2 + p_2 +\frac{3 s}{2}-\frac{C_1}{2},\\
   z =&2 q_1^3 + 2 p_2 q_1 - (C_1 -3 s) q_1 + q_2,\\
   w =&6 q_1^4 + 8 p_2 q_1^2 + 2 q_1 q_2 - 4 (C_1 - 3 s) q_1^2 - p_2^2 + p_1+ \frac{t}{2} -C_1 s + \frac{3 s^2}{2} + \frac{C_1^2}{2} - \frac{C_2}{2}
   \end{aligned}
  \right. 
\end{equation}
take the Hamiltonian system \eqref{eq:3} to the system
\begin{equation}\label{eq:1111112}
  \left\{
  \begin{aligned}
   dx =&y dt+\left(w - 6 y x^2 + y (C_1 - 3 s) \right)ds,\\
   dy =&z dt+h_1(x,y,z,w) ds,\\
   dz =&w dt+h_2(x,y,z,w) ds,\\
   dw =&(10 y^2 x + 10 x^2 z + 2 (C_1 - 3 s) x^3 - 6 x^5 + 
 \frac{1}{2} (-3 s^2 + 2 C_1 s +2 t + C_1^2 - 2 C_2) x\\
 & - (C_1 - 3 s) z-\alpha_0 + \frac{1}{2} )dt+h_3(x,y,z,w) ds,
   \end{aligned}
  \right. 
\end{equation}
where $h_i(x,y,z,w) \in {\mathbb C}(t,s)[x,y,z,w] \hspace{0.2cm} (i=1,2,3)$.
\end{theorem}
Setting $u:=x$, we see that
\begin{equation}
\frac{\partial u}{\partial t}=y, \quad  \frac{\partial^2 u}{\partial t^2}=z, \quad \frac{\partial^3 u}{\partial t^3}=w,
\end{equation}
and
\begin{equation}\label{eq:1111113}
  \left\{
  \begin{aligned}
   \frac{\partial^4 u}{\partial t^4} =&10u \left(\frac{\partial u}{\partial t} \right)^2+ 10 u^2 \frac{\partial^2 u}{\partial t^2} + 2 (C_1 - 3 s) u^3 - 6 u^5\\
&+ \frac{1}{2} (-3 s^2 + 2 C_1 s +2 t + C_1^2 - 2 C_2)u- (C_1 -3 s) \frac{\partial^2 u}{\partial t^2} -\alpha_0 + \frac{1}{2},\\
   \frac{\partial u}{\partial s} =&\frac{\partial^3 u}{\partial t^3}-6u^2\frac{\partial u}{\partial t}+(C_1 - 3 s) \frac{\partial u}{\partial t}.
   \end{aligned}
  \right. 
\end{equation}
Setting $s=0$ and $C_1=C_2=0$, the first equation in \eqref{eq:1111113} just coincides with the second member of the second Painlev\'e hierarchy $P_{II}^{(2)}$. 

\noindent
The second equation can be considered as a modified type of the mKdV equation.

\section{Autonomous version of the system  \eqref{eq:3} and mKdV5 equation}

In this section, we find an autonomous version of the system  \eqref{eq:3} given by
\begin{equation}\label{eq:10}
  \left\{
  \begin{aligned}
   dq_1 =&\frac{\partial K_1}{\partial p_1}dt+\frac{\partial K_2}{\partial p_1}ds,\\
   dp_1 =&-\frac{\partial K_1}{\partial q_1}dt-\frac{\partial K_2}{\partial q_1}ds,\\
   dq_2 =&\frac{\partial K_1}{\partial p_2}dt+\frac{\partial K_2}{\partial p_2}ds,\\
   dp_2 =&-\frac{\partial K_1}{\partial q_2}dt-\frac{\partial K_2}{\partial q_2}ds
   \end{aligned}
  \right. 
\end{equation}
with the polynomial Hamiltonians
\begin{align}\label{eq:ph1}
\begin{split}
K_1=&q_1^2p_1+\alpha_0 q_1-\frac{q_2^2}{2}-p_2^3-\frac{3}{20}p_2+p_1 p_2,\\
K_2=&\frac{p_1^2}{2}-\frac{3}{20}p_1-\alpha_1 q_2-p_1p_2^2+2q_1^2 p_1p_2+2q_1p_1q_2+2\alpha_0 q_1p_2.
\end{split}
\end{align}

\begin{proposition}
The system \eqref{eq:10} satisfies the compatibility conditions$:$
\begin{equation}
\frac{\partial }{\partial s} \frac{\partial q_1}{\partial t}=\frac{\partial }{\partial t} \frac{\partial q_1}{\partial s}, \quad \frac{\partial }{\partial s} \frac{\partial p_1}{\partial t}=\frac{\partial }{\partial t} \frac{\partial p_1}{\partial s}, \quad \frac{\partial }{\partial s} \frac{\partial q_2}{\partial t}=\frac{\partial }{\partial t} \frac{\partial q_2}{\partial s}, \quad \frac{\partial }{\partial s} \frac{\partial p_2}{\partial t}=\frac{\partial }{\partial t} \frac{\partial p_2}{\partial s}.
\end{equation}
\end{proposition}

\begin{proposition}
The system \eqref{eq:10} has $K_1$ and $K_2$ as its first integrals.
\end{proposition}

\begin{proposition}
Two Hamiltonians $K_1$ and $K_2$ satisfy
\begin{equation}
\{K_1,K_2\}=0,
\end{equation}
where
\begin{equation}
\{K_1,K_2\}=\frac{\partial K_1}{\partial p_1}\frac{\partial K_2}{\partial q_1}-\frac{\partial K_1}{\partial q_1}\frac{\partial K_2}{\partial p_1}+\frac{\partial K_1}{\partial p_2}\frac{\partial K_2}{\partial q_2}-\frac{\partial K_1}{\partial q_2}\frac{\partial K_2}{\partial p_2}.
\end{equation}
\end{proposition}
Here, $\{,\}$ denotes the poisson bracket such that $\{p_i,q_j\}={\delta}_{ij}$ (${\delta}_{ij}$:kronecker's delta).

\begin{proposition}
The system \eqref{eq:10} admits the extended affine Weyl group symmetry of type $A_1^{(1)}$ as the group of its B{\"a}cklund transformations, whose generators $s_0,s_1,{\pi}$ defined as follows$:$ with {\it the notation} $(*):=(q_1,p_1,q_2,p_2,t,s;\alpha_0,\alpha_1)$\rm{: \rm}
\begin{align*}
s_0:(*) \rightarrow &\left(q_1+\frac{\alpha_0}{f_0},p_1,q_2,p_2,t,s;-\alpha_0,\alpha_1+2\alpha_0 \right),\\
s_1:(*) \rightarrow &(q_1+\frac{\alpha_1}{f_1},p_1-\frac{4\alpha_1(q_2+2q_1p_2)}{f_1}+\frac{4\alpha_1^2(p_2+2q_1^2)}{f_1^2},\\
&q_2+\frac{4\alpha_1(q_1^2-p_2)}{f_1}+\frac{12\alpha_1^2 q_1}{f_1^2}+\frac{4\alpha_1^3}{f_1^3},p_2-\frac{4\alpha_1 q_1}{f_1}-\frac{2\alpha_1^2}{f_1^2},t,s;\alpha_0+2\alpha_1,-\alpha_1),\\
\pi:(*) \rightarrow &(-q_1,-f_1,-(q_2+4q_1p_2+4q_1^3),-(p_2+2q_1^2),t,s;\alpha_1,\alpha_0),
\end{align*}
where $f_0:=p_1$ and $f_1:=p_1+4q_1^2p_2-2p_2^2+4q_1q_2-\frac{3}{10}$.
\end{proposition}
Here, the parameters $\alpha_i$ satisfy the relation $\alpha_0+\alpha_1=0$.

\begin{proposition}
The system \eqref{eq:10} admits a rational solution$:$
\begin{equation}
(q_1,p_1,q_2,p_2;\alpha_0,\alpha_1)=\left(0,\frac{3}{20},0,0;0,0 \right).
\end{equation}
\end{proposition}

\begin{proposition}
Let us consider a polynomial Hamiltonian system with Hamiltonian $K \in {\mathbb C}[q_1,p_1,q_2,p_2]$. We assume that

$(C1)$ $deg(K)=4$ with respect to $q_1,p_1,q_2,p_2$.

$(C2)$ This system becomes again a polynomial Hamiltonian system in each coordinate $R_i \ (i=0,1)${\rm : \rm}
\begin{align*}
\begin{split}
R_0:(x_0,y_0,z_0,w_0)=&\left(\frac{1}{q_1},-(q_1f_0+\alpha_0)q_1,q_2,p_2 \right),\\
R_1:(x_1,y_1,z_1,w_1)=&\left(\frac{1}{q_1},-(q_1f_1+\alpha_1)q_1,q_2+4q_1p_2+4q_1^3,p_2+2q_1^2 \right),
\end{split}
\end{align*}
where $f_0:=p_1$ and $f_1:=p_1+4q_1^2p_2-2p_2^2+4q_1q_2-\frac{3}{10}$.

$(C3)$ $\alpha_0+\alpha_1=0$.

\noindent
Then such a system coincides with the Hamiltonian system \eqref{eq:10} with the polynomial Hamiltonians $K_1,K_2$.
\end{proposition}
We note that the conditions $(C2)$ should be read that
\begin{align*}
\begin{split}
&R_0(K), \quad R_1(K)
\end{split}
\end{align*}
are polynomials with respect to $x_i,y_i,z_i,w_i$.

Let us consider the relation between the polynomial Hamiltonian system \eqref{eq:10} and mKdV5 equation. At first, we can make the birational transformations between the polynomial Hamiltonian system \eqref{eq:10} and mKdV5 equation.
\begin{theorem}
The birational transformations
\begin{equation}\label{eq:11}
  \left\{
  \begin{aligned}
   x =&-q_1,\\
   y =&-(p_2+q_1^2),\\
   z =&-(q_2+2q_1p_2+2q_1^3),\\
   w =&-\left(p_1+2q_1q_2-p_2^2+8q_1^2p_2+6q_1^4-\frac{3}{20} \right)
   \end{aligned}
  \right. 
\end{equation}
take the Hamiltonian system \eqref{eq:10} to the system
\begin{equation}\label{eq:12}
  \left\{
  \begin{aligned}
   dx =&y dt+(w-6x^2y)ds,\\
   dy =&z dt+\left(4x^2z-2xy^2-6x^5-\frac{3}{10}x+\alpha_0 \right)ds,\\
   dz =&w dt+\left(-30x^4y-2y^3+4xyz+4x^2w-\frac{3}{10}y \right)ds,\\
   dw =&\left(10xy^2+10x^2z-6x^5-\frac{3}{10}x+\alpha_0 \right)dt\\
&+(-24x^7+10x^4z-80x^3y^2+10x^4z+12xyw\\
&-2y^2z+4xz^2-\frac{6}{5}x^3+4\alpha_0 x^2-\frac{3}{10}z)ds.
   \end{aligned}
  \right. 
\end{equation}
\end{theorem}
Setting $u:=x$, we see that
\begin{equation}
\frac{\partial u}{\partial t}=y, \quad  \frac{\partial^2 u}{\partial t^2}=z, \quad \frac{\partial^3 u}{\partial t^3}=w,
\end{equation}
and
\begin{equation}\label{eq:13}
  \left\{
  \begin{aligned}
   \frac{\partial^4 u}{\partial t^4} =&10u\left(\frac{\partial u}{\partial t} \right)^2+10u^2\frac{\partial^2 u}{\partial t^2}-6u^5-\frac{3}{10}u+\alpha_0,\\
   \frac{\partial u}{\partial s} =&\frac{\partial^3 u}{\partial t^3}-6u^2\frac{\partial u}{\partial t}.
   \end{aligned}
  \right. 
\end{equation}
The first equation in \eqref{eq:13} coincides with an autonomous version of $P_{II}^{(2)}$, and the second equation just coincides with the mKdV equation. The independent variables $t$ and $s$ in the second equation coincides with the ones of the polynomial Hamiltonian system \eqref{eq:10}.

We see that the first equation in \eqref{eq:13} is a compatible vector field for the mKdV equation.

Making a partial derivation in the variable $t$ for the first equation of \eqref{eq:13}, we obtain
\begin{equation}\label{eq:14}
  \left\{
  \begin{aligned}
   \frac{\partial^5 u}{\partial t^5} =&-30u^4\frac{\partial u}{\partial t}+10\left(\frac{\partial u}{\partial t} \right)^3+40u\frac{\partial u}{\partial t}\frac{\partial^2 u}{\partial t^2}+10u^2\frac{\partial^3 u}{\partial t^3}-\frac{3}{10}\frac{\partial u}{\partial t},\\
   0=&\frac{\partial u}{\partial s}-\frac{\partial^3 u}{\partial t^3}+6u^2\frac{\partial u}{\partial t}.
   \end{aligned}
  \right. 
\end{equation}

Adding each system in \eqref{eq:14}, we can obtain mKdV5 equation
\begin{equation}
   \frac{\partial^5 u}{\partial t^5} =10\left(u^2-\frac{1}{10} \right) \frac{\partial^3 u}{\partial t^3}+40u\frac{\partial u}{\partial t}\frac{\partial^2 u}{\partial t^2}+10\left(\frac{\partial u}{\partial t} \right)^3-30\left(u^2-\frac{1}{10} \right)^2 \frac{\partial u}{\partial t}+\frac{\partial u}{\partial s}.
\end{equation}

Now, let us solve the polynomial Hamiltonian system \eqref{eq:10} in the variable $p_2$.
\begin{theorem}
The birational transformations
\begin{equation}\label{eq:a1}
  \left\{
  \begin{aligned}
   x_1 =&-2q_1 p_1-6q_2 p_2-\alpha_0,\\
   y_1 =&p_1-3p_2^2-\frac{3}{20},\\
   z_1 =&q_2,\\
   w_1 =&p_2
   \end{aligned}
  \right. 
\end{equation}
take the Hamiltonian system \eqref{eq:10} to the system
\begin{equation}\label{eq:a2}
  \left\{
  \begin{aligned}
   dw_1 =&z_1 dt+(x_1+6z_1 w_1)ds,\\
   dz_1 =&y_1 dt+f_1(x_1,y_1,z_1,w_1)ds,\\
   dy_1 =&x_1 dt+f_2(x_1,y_1,z_1,w_1)ds,\\
   dx_1 =&(-\frac{1}{10(3+60w_1^2+20y_1)}(9w_1+360w_1^3+3600w_1^5-100x_1^2+300y_1 w_1+6000y w_1^3\\
&+1600y_1^2 w_1-1200 x_1 z_1 w_1+180z_1^2+1200y_1 z_1^2+100 \alpha_0^2)dt+f_3(x_1,y_1,z_1,w_1)ds,
   \end{aligned}
  \right. 
\end{equation}
where $f_i(x_1,y_1,z_1,w_1) \in {\mathbb C}(x_1,y_1,z_1,w_1)$.
\end{theorem}
Setting $v:=w_1$, we see that
\begin{equation}
\frac{\partial v}{\partial t}=z_1, \quad  \frac{\partial^2 v}{\partial t^2}=y_1, \quad \frac{\partial^3 v}{\partial t^3}=x_1,
\end{equation}
and
\begin{equation}\label{eq:a3}
  \left\{
  \begin{aligned}
   &2 \left(\frac{\partial^2 v}{\partial t^2}+3v^2+\frac{3}{20} \right) \left(\frac{\partial^4 v}{\partial t^4}+6\left(\frac{\partial v}{\partial t} \right)^2+6v\frac{\partial^2 v}{\partial t^2} \right)\\
&-\left(\frac{\partial^3 v}{\partial t^3}+6v\frac{\partial v}{\partial t} \right)^2+4v \left(\frac{\partial^2 v}{\partial t^2}+3v^2+\frac{3}{20} \right)^2 =-\alpha_0^2,\\
   &\frac{\partial v}{\partial s} =\frac{\partial^3 v}{\partial t^3}+6v\frac{\partial v}{\partial t}.
   \end{aligned}
  \right. 
\end{equation}
The first equation in \eqref{eq:a3} coincides with an autonomous version of the second member of the Painlev\'e $P34$ hierarchy, and the second equation just coincides with the KdV equation. We remark that the independent variables $t$ and $s$ in the second equation coincides with the ones of the polynomial Hamiltonian system \eqref{eq:10}.

\begin{theorem}
The rational transformations
\begin{equation}\label{eq:a4}
  \left\{
  \begin{aligned}
   x_1 =&\frac{1}{10}(3x-100xw^2+60x^5-20xy-60zw-100x^2 z-10\alpha_0),\\
   y_1 =&-2w^2-y-2xz,\\
   z_1 =&-2xw-z,\\
   w_1 =&-x^2-w
   \end{aligned}
  \right. 
\end{equation}
take the system \eqref{eq:13} to the system \eqref{eq:a3}.
\end{theorem}
The fourth transformation in \eqref{eq:a4} just coincides with the Miura transformation:
\begin{equation}\label{eq:a5}
   v =-\left(\frac{\partial u}{\partial t}+u^2 \right),
\end{equation}
where $u:=x, w_1:=v$.

Next, let us solve the polynomial Hamiltonian system \eqref{eq:10} in the variable $p_1$.
\begin{theorem}
The birational transformations
\begin{equation}\label{eq:b1}
  \left\{
  \begin{aligned}
   x_2 =&-2q_1 p_1-\alpha_0,\\
   y_2 =&p_1,\\
   z_2 =&-2(-4q_1p_1p_2+p_1q_2-2\alpha_0 p_2),\\
   w_2 =&-2(p_1p_2-q_1^2 p_1-\alpha_0 q_1)
   \end{aligned}
  \right. 
\end{equation}
take the Hamiltonian system \eqref{eq:10} to the system
\begin{equation}\label{eq:b2}
  \left\{
  \begin{aligned}
   dy_2 =&x_2 dt+\left(\frac{3x_2(x_2-\alpha_0)(x_2+\alpha_0)}{2y_2^2}-\frac{3}{y_2}x_2w_2+z_2 \right)ds,\\
   dx_2 =&w_2 dt+f_1(x_2,y_2,z_2,w_2)ds,\\
   dw_2 =&z_2 dt+f_2(x_2,y_2,z_2,w_2)ds,\\
   dz_2 =&(\frac{3(3x_2-\alpha_0)(3x_2+\alpha_0)(x_2+\alpha_0)(x_2-\alpha_0)}{8y_2^3}-\frac{w_2(17x_2^2-5\alpha_0^2)}{2y_2^2}\\
&+\frac{7w_2^2+6x_2z_2}{2y_2}-2y_2^2+\frac{3}{10}y_2)dt+f_3(x_2,y_2,z_2,w_2)ds,
   \end{aligned}
  \right. 
\end{equation}
where $f_i(x_2,y_2,z_2,w_2) \in {\mathbb C}(x_2,y_2,z_2,w_2)$.
\end{theorem}
Setting $U:=y_2$, we see that
\begin{equation}
\frac{\partial U}{\partial t}=x_2, \quad  \frac{\partial^2 U}{\partial t^2}=w_2, \quad \frac{\partial^3 U}{\partial t^3}=z_2,
\end{equation}
and
\begin{equation}\label{eq:b3}
  \left\{
  \begin{aligned}
   \frac{\partial^4 U}{\partial t^4}=&\frac{3(3\frac{\partial U}{\partial t}-\alpha_0)(3\frac{\partial U}{\partial t}+\alpha_0)(\frac{\partial U}{\partial t}+\alpha_0)(\frac{\partial U}{\partial t}-\alpha_0)}{8U^3}-\frac{\frac{\partial^2 U}{\partial t^2}(17 \left(\frac{\partial U}{\partial t} \right)^2-5\alpha_0^2)}{2U^2}\\
&+\frac{7\left(\frac{\partial^2 U}{\partial t^2} \right)^2+6\frac{\partial U}{\partial t}\frac{\partial^3 U}{\partial t^3}}{2U}-2U^2+\frac{3}{10}U,\\
   \frac{\partial U}{\partial s} =&\frac{3\frac{\partial U}{\partial t}(\frac{\partial U}{\partial t}-\alpha_0)(\frac{\partial U}{\partial t}+\alpha_0)}{2U^2}-\frac{3}{U}\frac{\partial U}{\partial t}\frac{\partial^2 U}{\partial t^2}+\frac{\partial^3 U}{\partial t^3}.
   \end{aligned}
  \right. 
\end{equation}
Setting
\begin{equation}
deg(U)=1, \quad deg\left(\frac{\partial U}{\partial t} \right)=2, \quad  deg\left(\frac{\partial^2 U}{\partial t^2} \right)=3, \quad deg \left(\frac{\partial^3 U}{\partial t^3} \right)=4,
\end{equation}
\begin{equation}
deg\left(\frac{\partial U}{\partial s} \right)=4, \quad deg(\alpha_0)=2,
\end{equation}
the second equation in \eqref{eq:b3} can be considered as homogeneous equation of degree 4.

For the second equation in \eqref{eq:b3}, we see that this system admits travelling wave solutions $u(t,s)=U(t+cs)$, where $U(T) \ (T:=t+cs, \ c \in {\mathbb C})$ satisfies the equation
\begin{equation}
c\frac{d U}{d T} =\frac{3\frac{d U}{d T} \left(\frac{d U}{d T}-\alpha_0 \right) \left(\frac{d U}{d T}+\alpha_0 \right)}{2U^2}-\frac{3}{U}\frac{d U}{d T}\frac{d^2 U}{d T^2}+\frac{d^3 U}{\partial T^3}.
\end{equation}
After integrating once, this equation becomes
\begin{equation}
\frac{d^2 U}{d T^2}=\frac{3}{2U}\left(\frac{dU}{dT} \right)^2-\frac{3\alpha_0^2}{2U}+c U.
\end{equation}
Setting $x:=U,y:=\frac{dU}{dT}$, the birational transformations
\begin{equation}
  \left\{
  \begin{aligned}
   x_1 =&\frac{1}{x},\\
   y_1 =&\frac{y+\alpha_0}{x}
   \end{aligned}
  \right. 
\end{equation}
take the system
\begin{equation}
  \left\{
  \begin{aligned}
   \frac{dx}{dT} =&y ,\\
   \frac{dy}{dT} =&\frac{3y^2}{2x}-\frac{3\alpha_0^2}{2x}+c x
   \end{aligned}
  \right. 
\end{equation}
to the Hamiltonian system
\begin{equation}
  \left\{
  \begin{aligned}
   \frac{dx_1}{dT} =&\frac{\partial K}{\partial y_1}=-x_1y_1+\alpha_0 x_1^2 ,\\
   \frac{dy_1}{dT} =&-\frac{\partial K}{\partial x_1}=\frac{1}{2}y_1^2-2\alpha_0x_1y_1+c
   \end{aligned}
  \right. 
\end{equation}
with the polynomial Hamiltonian $K$
\begin{equation}
K:=-\frac{1}{2}x_1y_1^2+\alpha_0 x_1^2 y_1-cx_1.
\end{equation}
Elimination of $y_1$ from this system gives the second-order ordinary differential equation for $x_1$
\begin{equation}
\frac{d^2x_1}{dT^2} =\frac{1}{2x_1} \left(\frac{dx_1}{dT} \right)^2+\frac{3}{2} \alpha_0^2 x_1^3-c x_1.
\end{equation}
It is well-known that this equation is an autonomous version of the fourth Painlev\'e equation.

Next, let us solve the polynomial Hamiltonian system \eqref{eq:10} in the variables $p_1$ and $p_2$.
\begin{theorem}
The birational transformations
\begin{equation}\label{eq:c1}
  \left\{
  \begin{aligned}
   x_3 =&-2q_1 p_1-\alpha_0,\\
   y_3 =&p_1,\\
   z_3 =&q_2,\\
   w_3 =&p_2
   \end{aligned}
  \right. 
\end{equation}
take the Hamiltonian system \eqref{eq:10} to the system
\begin{equation}\label{eq:c2}
  \left\{
  \begin{aligned}
   dx_3 =&\left(\frac{(x_3+\alpha_0)(x_3-\alpha_0)}{2y_3}-2y_3w_3 \right) dt\\
&+\left(\frac{w_2(x_3+\alpha_0)(x_3-\alpha_0)}{2y_3}+\frac{3}{10}y_3+2y_3 w_3^2-2y_3^2 \right)ds,\\
   dy_3 =&x_3 dt+(2x_3 w_3-2y_3z_3)ds,\\
   dz_3 =&\left(-3w_3^2+y_3-\frac{3}{20} \right) dt+\left(\frac{(x_3+\alpha_0)(x_3-\alpha_0)}{2y_3}-2y_3w_3 \right)ds,\\
   dw_3 =&z_3dt+x_3ds.
   \end{aligned}
  \right. 
\end{equation}
\end{theorem}
Setting $u_1:=y_3$ and $v_1:=w_3$, we see that
\begin{equation}
\frac{\partial u_1}{\partial t}=x_3, \quad  \frac{\partial v_1}{\partial t}=z_3, \quad  \frac{\partial v_1}{\partial s}=x_3=\frac{\partial u_1}{\partial t}
\end{equation}
and
\begin{equation}\label{eq:c3}
  \left\{
  \begin{aligned}
   \frac{\partial^2 u_1}{\partial t^2}=&\frac{(\frac{\partial u_1}{\partial t}+\alpha_0)(\frac{\partial u_1}{\partial t}-\alpha_0)}{2u_1}-2u_1v_1,\\
   \frac{\partial^2 v_1}{\partial t^2}=&-3v_1^2+u_1-\frac{3}{20}
   \end{aligned}
  \right. 
\end{equation}
and
\begin{equation}\label{eq:c4}
 \frac{\partial u_1}{\partial s}=2\frac{\partial u_1}{\partial t} v_1-2u_1 \frac{\partial v_1}{\partial t}.
\end{equation}

Finally, let us solve the polynomial Hamiltonian system \eqref{eq:10} in the variables $q_1$ and $q_2$.
\begin{theorem}
The birational transformations
\begin{equation}\label{eq:g1}
  \left\{
  \begin{aligned}
   x_4 =&q_1,\\
   y_4 =&p_1-3p_2^2-\frac{3}{20},\\
   z_4 =&q_2,\\
   w_4 =&p_2+q_1^2
   \end{aligned}
  \right. 
\end{equation}
take the Hamiltonian system \eqref{eq:10} to the system
\begin{equation}\label{eq:g2}
  \left\{
  \begin{aligned}
   dx_4 =&w_4dt+\left(2w_4^2-2x_4^2 w_4+y_4+2x_4z_4 \right)ds,\\
   dw_4 =&(z_4+2x_4w_4) dt+h_1(x_4,y_4,z_4,w_4)ds,\\
   dz_4 =&y_4 dt+\left(-\frac{3}{10}w_4-6w_4^3+\frac{3}{5}x_4^2+24x_4^2 w_4^2-30x_4^4 w_4+12x_4^6-2y_4w_4+4x_4^2 y_4+2\alpha_0 x_4 \right)ds,\\
   dy_4 =&\left(-\frac{3}{10}x_4-6x_4w_4^2+12x_4^3 w_4-6x_4^5-2x_4y_4-6z_4w_4+6x_4^2 z_4-\alpha_0 \right)dt\\
&+h_2(x_4,y_4,z_4,w_4)ds,
   \end{aligned}
  \right. 
\end{equation}
where $h_i(x_4,y_4,z_4,w_4) \in {\mathbb C}[x_4,y_4,z_4,w_4]$.
\end{theorem}
Setting $u_2:=x_4$ and $v_2:=z_4$, we see that
\begin{equation}
\frac{\partial u_2}{\partial t}=w_4, \quad  \frac{\partial v_2}{\partial t}=y_4,
\end{equation}
and
\begin{equation}\label{eq:g3}
  \left\{
  \begin{aligned}
   \frac{\partial^2 u_2}{\partial t^2}=&v_2+2u_2\frac{\partial u_2}{\partial t},\\
   \frac{\partial^2 v_2}{\partial t^2}=&-\frac{3}{10}u_2-6u_2 \left(\frac{\partial u_2}{\partial t} \right)^2+12u_2^3 \frac{\partial u_2}{\partial t}-6u_2^5-2u_2\frac{\partial v_2}{\partial t}-6v_2 \frac{\partial u_2}{\partial t}+6u_2^2 v_2-\alpha_0
   \end{aligned}
  \right. 
\end{equation}
and
\begin{equation}\label{eq:g4}
  \left\{
  \begin{aligned}
   \frac{\partial u_2}{\partial s}=&2\left(\frac{\partial u_2}{\partial t} \right)^2-2u_2^2 \frac{\partial u_2}{\partial t}+\frac{\partial v_2}{\partial t}+2u_2v_2,\\
   \frac{\partial v_2}{\partial s}=&-\frac{3}{10} \frac{\partial u_2}{\partial t}-6\left(\frac{\partial u_2}{\partial t} \right)^3+\frac{3}{5}u_2^2+24u_2^2 \left(\frac{\partial u_2}{\partial t} \right)^2-30u_2^4 \frac{\partial u_2}{\partial t}+12u_2^6-2\frac{\partial v_2}{\partial t} \frac{\partial u_2}{\partial t}\\
&+4u_2^2 \frac{\partial v_2}{\partial t}+2\alpha_0 u_2.
   \end{aligned}
  \right. 
\end{equation}
Setting
\begin{equation}
deg(u_2)=1, \quad deg(v_2)=3, \quad deg\left(\frac{\partial u_2}{\partial t} \right)=2, \quad deg\left(\frac{\partial v_2}{\partial t} \right)=4, \quad  deg\left(\frac{\partial u_2}{\partial s} \right)=4,
\end{equation}
the first equation in \eqref{eq:g4} can be considered as homogeneous equation of degree 4.

\section{MkdV equation and Lax pair}
In this section, we study the following differential equations in three variables $T,t,s$
\begin{equation}\label{lx1}
  \left\{
  \begin{aligned}
   {\partial}_{t}(A(T,t,s))&=[B_1(T,t,s),A(T,t,s)](=B_1(T,t,s)A(T,t,s)-A(T,t,s)B_1(T,t,s)),\\
   {\partial}_{s}(A(T,t,s))&=[B_3(T,t,s),A(T,t,s)](=B_3(T,t,s)A(T,t,s)-A(T,t,s)B_3(T,t,s)),
   \end{aligned}
  \right. 
\end{equation}
where the matrices $A,B_1$ and $B_3$ (cf. \cite{N1,N2,N3,N4}, see Section 2) are given by
\begin{align}
\begin{split}
&-A(T,t,s)=\begin{pmatrix}
\varepsilon_1 & y\\
0 & \varepsilon_2 
\end{pmatrix}+\begin{pmatrix}
2z & 4q\\
-x & -2z 
\end{pmatrix}T+\begin{pmatrix}
8p & -8\\
-4w & -8p
\end{pmatrix}T^2+\begin{pmatrix}
0 & 0\\
-8 & 0 
\end{pmatrix}T^3,\\
&B_1(T,t,s)=\begin{pmatrix}
p & -1\\
0 & -p 
\end{pmatrix}+\begin{pmatrix}
0 & 0\\
-1 & 0
\end{pmatrix}T,\\
&B_3(T,t,s)=\begin{pmatrix}
\frac{z}{4} & \frac{q}{2}\\
0 & -\frac{z}{4} 
\end{pmatrix}+\begin{pmatrix}
p & -1\\
-\frac{w}{2} & -p 
\end{pmatrix}T+\begin{pmatrix}
0 & 0\\
-1 & 0
\end{pmatrix}T^2,
\end{split}
\end{align}
where $x,y,z,w,q,p$ denote unknown complex variables in $t,s$, and $\varepsilon_i$ are constant complex parameters. These Lax pairs are well-known as the ones of the soliton equations.

Considering the autonomous limit $\delta=0$ for the following equations

\begin{align}
\begin{split}
&{\partial}_{t}(A(T,t,s))-\delta T{\partial}_{T}(B_1(T,t,s))+[A(T,t,s),B_1(T,t,s)]=0,\\
&{\partial}_{s}(A(T,t,s))-\delta T{\partial}_{T}(B_3(T,t,s))+[A(T,t,s),B_3(T,t,s)]=0,\\
&\delta T{\partial}_{T} \vec{u}=A(T,t,s)\vec{u}, \quad {\partial}_{t} \vec{u}=B_1(T,t,s)\vec{u}, \quad {\partial}_{s} \vec{u}=B_3(T,t,s)\vec{u}, \quad \vec{u}=\begin{pmatrix}
u_1\\
u_2 
\end{pmatrix},
\end{split}
\end{align}
we can obtain the above system \eqref{lx1}.

By solving the system \eqref{lx1}, we will obtain the Hamiltonian system \eqref{eq:10} given in Section 5. Here, $\alpha_0,\alpha_1$ are complex parameters satisfying the relation:
\begin{equation}
  \left\{
  \begin{aligned}
  \alpha_0  &=-\varepsilon_1+\varepsilon_2,\\
  \alpha_1  &=\varepsilon_1-\varepsilon_2,\\
  \alpha_0+\alpha_1&=0.
   \end{aligned}
  \right. 
\end{equation}

The latter part of this section is devoted to showing the above results.

At first, we study a 1-parameter family of total differential system in two variables $t,s$
\begin{equation}\label{eq:aa1}
  \left\{
  \begin{aligned}
   dx =&(-2xp-\alpha_0)dt+(-2xz-2\alpha_0 w)ds,\\
   dy =&(2yp+\alpha_1)dt+(2yz-2\alpha_1 q)ds,\\
   dz =&\left(\frac{x+y}{2} \right)dt+(-xq+yw-2\delta p)ds,\\
   dw =&(z-2wp)dt+(-2xp+\alpha_1-2\delta)ds,\\
   dq =&(z+2qp)dt+(2yp+\alpha_1+\delta)ds,\\
   dp =&\left(\frac{w+q}{2} \right)dt+\left(\frac{x+y}{2} \right)ds,
   \end{aligned}
  \right.
\end{equation}
where the constant parameters $\alpha_i$ satisfy the relation:
\begin{equation}
\alpha_0+\alpha_1=\delta \quad (\delta \in {\mathbb C}).
\end{equation}
This system can be considered as a modified version of the system \eqref{eq:1} with the constant parameter $\delta$.

We similarly show that the system \eqref{eq:aa1} can be obtained by the compatibility conditions for the second-order linear differential equations in three variables $T,t,s$
\begin{equation}
\delta T{\partial}_{T} \vec{u}=A(T,t,s)\vec{u}, \quad {\partial}_{t} \vec{u}=B_1(T,t,s)\vec{u}, \quad \vec{u}=\begin{pmatrix}
u_1\\
u_2 
\end{pmatrix},
\end{equation}
and
\begin{equation}
\delta T{\partial}_{T} \vec{u}=A(T,t,s)\vec{u}, \quad {\partial}_{s} \vec{u}=B_3(T,t,s)\vec{u}.
\end{equation}

Now, let us consider the case $\delta=0$. In this case, we consider the following system
\begin{equation}\label{eq:a111}
  \left\{
  \begin{aligned}
   dx =&(-2xp-\alpha_0)dt+(-2xz-2\alpha_0 w)ds,\\
   dy =&(2yp+\alpha_1)dt+(2yz-2\alpha_1 q)ds,\\
   dz =&\left(\frac{x+y}{2} \right)dt+(-xq+yw)ds,\\
   dw =&(z-2wp)dt+(-2xp+\alpha_1)ds,\\
   dq =&(z+2qp)dt+(2yp+\alpha_1)ds,\\
   dp =&\left(\frac{w+q}{2} \right)dt+\left(\frac{x+y}{2} \right)ds.
   \end{aligned}
  \right.
\end{equation}
We easily see that each equation
\begin{equation}
\frac{\partial z}{\partial s}=-xq+yw-2\delta p, \quad \frac{\partial w}{\partial s}=-2xp+\alpha_1-2\delta, \quad  \frac{\partial q}{\partial s}=2yp+\alpha_1+\delta
\end{equation}
is changed into
\begin{equation}
\frac{\partial z}{\partial s}=-xq+yw \quad \frac{\partial w}{\partial s}=-2xp+\alpha_1, \quad  \frac{\partial q}{\partial s}=2yp+\alpha_1.
\end{equation}
We see that the system \eqref{eq:a111} can be obtained by solving the following Lax equations:
\begin{equation}\label{Lax1}
{\partial}_{t}(A(T,t,s))=[B_1(T,t,s),A(T,t,s)]
\end{equation}
and
\begin{equation}\label{Lax2}
{\partial}_{s}(A(T,t,s))=[B_3(T,t,s),A(T,t,s)].
\end{equation}
The Lax pairs \eqref{Lax1} and \eqref{Lax2} are well-known as the ones of the soliton equations.

Next, let us show that the system \eqref{eq:a111} is equivalent to the polynomial Hamiltonian system \eqref{eq:10} in two variables $t,s$.

We see that this system \eqref{eq:a111} has its first integrals:
\begin{equation}\label{fi}
  \left\{
  \begin{aligned}
    w-q+2p^2=&0,\\
    4zp-2wq+x-y=&\frac{3}{10},
   \end{aligned}
  \right.
\end{equation}
where we select the integral constants $C_1,C_2$ as
\begin{align}
\begin{split}
C_1=0, \quad C_2=\frac{3}{10}.
\end{split}
\end{align}
We similarly show that the birational transformations
\begin{equation}
  \left\{
  \begin{aligned}
   q_1 =&p,\\
   p_1 =&x,\\
   q_2 =&z-2wp,\\
   p_2 =&w
   \end{aligned}
  \right. 
\end{equation}
take the system \eqref{eq:a111} to the Hamiltonian system \eqref{eq:10} with the polynomial Hamiltonians \eqref{eq:ph1} given in Section 5. We remark that the independent variables $t$ and $s$ in the system \eqref{eq:10} coincide with the variables of the Lax pairs \eqref{Lax1} and \eqref{Lax2}.

The relations between $x,y,z,w,q,p$ and $q_1,p_1,q_2,p_2$ are given by
\begin{equation}\label{relation1}
  \left\{
  \begin{aligned}
   x =&p_1,\\
   y =&4q_1^2 p_2-2p_2^2+4q_1q_2+p_1-\frac{3}{10},\\
   z =&q_2+2q_1 p_2,\\
   w =&p_2,\\
   q =&2q_1^2+p_2,\\
   p =&q_1.
   \end{aligned}
  \right.
\end{equation}

Here, let us calculate the determinant of the matrix $A$:
\begin{align}
\begin{split}
det(A)=&-64 T^5 - 32 T^4 (2 p^2 - q + w) -
 8 T^3 (-2 q w + x - y + 4 p z)\\
&+4 T^2 (q x + w y - z^2 - 2(\varepsilon_1 - \varepsilon_2)p) + 
 T (x y - 2(\varepsilon_1 -\varepsilon_2)z) + \varepsilon_1 \varepsilon_2 .
\end{split}
\end{align}
Next, let us calculate the characteristic polynomial of the matrix $A$:
\begin{align}
\begin{split}
\lambda^2+det(A)=0.
\end{split}
\end{align}
If the eigenvalues of the matrix $A$ are independent in the variables $t$ and $s$, we obtain the following conditions:
\begin{align}\label{cod1}
\begin{split}
2 p^2 - q + w=&K_3,\\
-2 q w + x - y + 4 p z=&K_4,
\end{split}
\end{align}
and
\begin{align}\label{cod2}
\begin{split}
x y - 2(\varepsilon_1 - \varepsilon_2)z=&K_1,\\
q x + w y - z^2 - 2(\varepsilon_1 - \varepsilon_2)p=&K_2,
\end{split}
\end{align}
where, $K_1,K_2,K_3,K_4 \in {\mathbb C}$. By using the relation \eqref{relation1}, we see that the first conditions \eqref{cod1} with $K_3=0,K_4=\frac{3}{10}$ are equivalent to the conditions \eqref{fi}, and the second conditions \eqref{cod2} just coincide with the polynomial Hamiltonians \eqref{eq:ph1} given in Section 5.

We see that the phase space of the partial differential system \eqref{eq:10} corresponds to the algebraic surface $\{(q_1,p_1,q_2,p_2)|f_1=f_2=0\}$, where $f_1$ and $f_2$ are explicitly given as follows:
\begin{equation}
  \left\{
  \begin{aligned}
   f_1(q_1,p_1,q_2,p_2)=&q_1^2p_1+\alpha_0 q_1-\frac{q_2^2}{2}-p_2^3-\frac{3p_2}{20}+p_1 p_2-K_1,\\
f_2(q_1,p_1,q_2,p_2)=&\frac{p_1^2}{2}-\frac{3p_1}{20}-\alpha_1 q_2-p_1p_2^2+2q_1^2 p_1p_2+2q_1p_1q_2+2\alpha_0 q_1p_2-K_2.
   \end{aligned}
  \right. 
\end{equation}
By solving the equation $f_1=0$ with respect to the variable $p_1$, we can obtain the hypersuface $\{(q_1,q_2,p_2) \in {\mathbb C}^3|F(q_1,q_2,p_2)=0\}$;
\begin{align}
\begin{split}
&F(q_1,q_2,p_2)=\\
&-1600 p_2^4 q_1^4 - 800 p_2^5 q_1^2 - 1600 p_2^3 q_1^3 q_2 - 240 p_2^2 q_1^4 - 
 800 q_1^3 q_2^3 + 400 p_2^6 - 400 p_2^2 q_1^2 q_2^2 - 1600 p_2^4 q_1 q_2\\
&-800 p_2 q_1^4 q_2^2 - 800 p_2 q_1 q_2^3 - 2400 \alpha_0 p_2^2 q_1^3 + 
 1600 \alpha_0 p_2 q_1^2 q_2 + 1600 \alpha_0 q_1^4 q_2 + 
 1600 \alpha_1 p_2 q_1^2 q_2\\
&+ 800 \alpha_1 q_1^4 q_2 - 240 p_2 q_1^3 q_2 - 
 1600 K_1 p_2 q_1^4 + 120 p_2^4 + 800 K_2 q_1^4 - 100 q_2^4 - 
 1600 \alpha_0 p_2^3 q_1\\
&- 800 K_1 p_2^2 q_1^2 - 1600 K_1 q_1^3 q_2 + 60 q_1^2 q_2^2 - 
 240 p_2^2 q_1 q_2 + 18 p_2 q_1^2 + 1600 K_2 p_2 q_1^2 + 800 \alpha_1 p_2^2 q_2\\
&+400 \alpha_0 q_1 q_2^2 - 1600 K_1 p_2 q_1 q_2 - 120 \alpha_0 q_1^3 - 
 400 \alpha_0^2 q_1^2 + 120 K_1 q_1^2 - 400 K_1 q_2^2 + 9 p_2^2\\
& + 800 K_2 p_2^2 +800 \alpha_0 K_1 q_1 - 400 K_1^2.
\end{split}
\end{align}
Here, we see that $deg(F)=8$ with respect to $q_1,q_2,p_2$.

\section{Autonomous version of $P_{II}^{(3)}$ and mKdV equation}

In this section, we find a one-parameter family of polynomial Hamiltonian system in two variables given by
\begin{equation}\label{eq:15}
  \left\{
  \begin{aligned}
   dq_1 =&\frac{\partial K_1}{\partial p_1}dt+\frac{\partial K_2}{\partial p_1}ds, \quad dp_1 =-\frac{\partial K_1}{\partial q_1}dt-\frac{\partial K_2}{\partial q_1}ds,\\
   dq_2 =&\frac{\partial K_1}{\partial p_2}dt+\frac{\partial K_2}{\partial p_2}ds, \quad dp_2 =-\frac{\partial K_1}{\partial q_2}dt-\frac{\partial K_2}{\partial q_2}ds,\\
   dq_3 =&\frac{\partial K_1}{\partial p_3}dt+\frac{\partial K_2}{\partial p_3}ds, \quad dp_3 =-\frac{\partial K_1}{\partial q_3}dt-\frac{\partial K_2}{\partial q_3}ds
   \end{aligned}
  \right. 
\end{equation}
with the polynomial Hamiltonians
\begin{align}
\begin{split}
K_1=&q_1^2p_1-\alpha_0 q_1-\frac{p_2^4}{2}-\frac{g}{2}p_2+\frac{1}{2}p_3^2+p_1 p_2-q_2q_3+p_2 q_3^2-2p_2^2 p_3,\\
K_2=&-\frac{1}{4}q_2^2-\frac{1}{4}gp_2^2+\frac{1}{2}q_3^2 p_3-\frac{1}{2} \alpha_0 q_3-\frac{1}{4}gp_3-\alpha_0 q_1p_2+\frac{1}{2}p_1p_3+q_1p_1q_3-p_2p_3^2\\
&+\frac{1}{2}p_2^2 q_3^2-p_2^3 p_3+q_1^2 p_1p_2 \quad (g \in {\mathbb C}).
\end{split}
\end{align}
Here, $\alpha_0$ and $\alpha_1$ are constant complex parameters.

\begin{proposition}
The system \eqref{eq:15} satisfies the compatibility conditions$:$
\begin{align}
\begin{split}
&\frac{\partial }{\partial s} \frac{\partial q_1}{\partial t}=\frac{\partial }{\partial t} \frac{\partial q_1}{\partial s}, \quad \frac{\partial }{\partial s} \frac{\partial p_1}{\partial t}=\frac{\partial }{\partial t} \frac{\partial p_1}{\partial s}, \quad \frac{\partial }{\partial s} \frac{\partial q_2}{\partial t}=\frac{\partial }{\partial t} \frac{\partial q_2}{\partial s},\\
&\frac{\partial }{\partial s} \frac{\partial p_2}{\partial t}=\frac{\partial }{\partial t} \frac{\partial p_2}{\partial s}, \quad \frac{\partial }{\partial s} \frac{\partial q_3}{\partial t}=\frac{\partial }{\partial t} \frac{\partial q_3}{\partial s}, \quad \frac{\partial }{\partial s} \frac{\partial p_3}{\partial t}=\frac{\partial }{\partial t} \frac{\partial p_3}{\partial s}.
\end{split}
\end{align}

\end{proposition}

\begin{proposition}
Two Hamiltonians $K_1$ and $K_2$ satisfy
\begin{equation}
\{K_1,K_2\}=0,
\end{equation}
where
\begin{equation}
\{K_1,K_2\}=\frac{\partial K_1}{\partial p_1}\frac{\partial K_2}{\partial q_1}-\frac{\partial K_1}{\partial q_1}\frac{\partial K_2}{\partial p_1}+\frac{\partial K_1}{\partial p_2}\frac{\partial K_2}{\partial q_2}-\frac{\partial K_1}{\partial q_2}\frac{\partial K_2}{\partial p_2}+\frac{\partial K_1}{\partial p_3}\frac{\partial K_2}{\partial q_3}-\frac{\partial K_1}{\partial q_3}\frac{\partial K_2}{\partial p_3}.
\end{equation}
\end{proposition}
Here, $\{,\}$ denotes the poisson bracket such that $\{p_i,q_j\}={\delta}_{ij}$ (${\delta}_{ij}$:kronecker's delta).

\begin{theorem}
The system \eqref{eq:15} has $K_1$ and $K_2$ as its first integrals.
\end{theorem}

\begin{theorem}
The system \eqref{eq:15} admits the extended affine Weyl group symmetry of type $A_1^{(1)}$ as the group of its B{\"a}cklund transformations, whose generators $s_0,s_1,{\pi}$ defined as follows$:$ with {\it the notation} $(*):=(q_1,p_1,q_2,p_2,q_3,p_3,t,s;\alpha_0,\alpha_1)$\rm{: \rm}
\begin{align*}
s_0:(*) \rightarrow &\left(q_1-\frac{\alpha_0}{f_0},p_1,q_2,p_2,q_3,p_3,t,s;-\alpha_0,\alpha_1+2\alpha_0 \right),\\
s_1:(*) \rightarrow &(q_1-\frac{\alpha_1}{f_1},p_1+\frac{4\alpha_1(q_2+2q_1p_3+2q_1 p_2^2)}{f_1}+\frac{4\alpha_1^2(p_3-p_2^2+4q_1q_3+4q_1^2p_2)}{f_1^2},\\
&q_2-\frac{4\alpha_1(p_3-2q_1^2 p_2)}{f_1}+\frac{8\alpha_1^2 (q_3+q_1p_2-2q_1^3)}{f_1^2}+\frac{8\alpha_1^3(p_2+2q_1^2)}{f_1^3},p_2+\frac{4\alpha_1 q_1}{f_1}-\frac{2\alpha_1^2}{f_1^2},\\
&q_3+\frac{4\alpha_1 (p_2-q_1^2)}{f_1}+\frac{12\alpha_1^2 q_1}{f_1^2}-\frac{4\alpha_1^3}{f_1^3},p_3+\frac{4\alpha_1 q_3}{f_1}+\frac{8\alpha_1^2 (p_2-q_1^2)}{f_1^2}+\frac{16\alpha_1^3 q_1}{f_1^3}-\frac{4\alpha_1^4}{f_1^4},\\
&t,s;\alpha_0+2\alpha_1,-\alpha_1),\\
\pi:(*) \rightarrow &(-q_1,-f_1,-(q_2+4q_1p_3+8q_1^2 q_3+8q_1^3 p_2),-(p_2+2q_1^2),-(q_3+4q_1p_2+4q_1^3),\\
&-(p_3+4q_1q_3+8q_1^2p_2+4q_1^4),t,s;\alpha_1,\alpha_0),
\end{align*}
where $f_0:=p_1$ and $f_1:=p_1+4q_1q_2+2q_3^2-4p_2p_3+4q_1^2 p_3+4q_1^2 p_2^2-g$.
\end{theorem}
Here, the parameters $\alpha_i$ satisfy the relation $\alpha_0+\alpha_1=0$.

\begin{proposition}
The system \eqref{eq:15} admits a rational solution$:$
\begin{equation}
(q_1,p_1,q_2,p_2,q_3,p_3;\alpha_0,\alpha_1)=\left(0,\frac{g}{2},0,0,0,0;0,0 \right).
\end{equation}
\end{proposition}

For the system \eqref{eq:15}, we find the holomorphy condition of this system. Thanks to this holomorphy condition, we can recover the Hamiltonian system \eqref{eq:15} with the polynomial Hamiltonians $K_1,K_2$.
\begin{theorem}
Let us consider a polynomial Hamiltonian system with Hamiltonian $K \in {\mathbb C}[q_1,p_1,q_2,p_2,q_3,p_3]$. We assume that

$(D1)$ $deg(K)=4$ with respect to $q_1,p_1,q_2,p_2,q_3,p_3$.

$(D2)$ This system becomes again a polynomial Hamiltonian system in each coordinate $R_i \ (i=0,1)${\rm : \rm}
\begin{align*}
\begin{split}
R_0:(x_0,y_0,z_0,w_0)=&\left(\frac{1}{q_1},-(q_1f_0-\alpha_0)q_1,q_2,p_2,q_3,p_3 \right),\\
R_1:(x_1,y_1,z_1,w_1)=&(\frac{1}{q_1},-(q_1f_1-\alpha_1)q_1,q_2+4q_1p_3+8q_1^2 q_3+8q_1^3 p_2,p_2+2q_1^2,\\
&q_3+4q_1p_2+4q_1^3,p_3+4q_1q_3+8q_1^2 p_2+4q_1^4 ),
\end{split}
\end{align*}
where $f_0:=p_1$ and $f_1:=p_1+4q_1q_2+2q_3^2-4p_2p_3+4q_1^2 p_3+4q_1^2 p_2^2-g$, and the parameters $\alpha_i$ satisfy the relation $\alpha_0+\alpha_1=0$. Then such a system coincides with the Hamiltonian system \eqref{eq:15} with the polynomial Hamiltonians $K_1,K_2$.
\end{theorem}
We note that the conditions $(D2)$ should be read that
\begin{align*}
\begin{split}
&R_0(K), \quad R_1(K)
\end{split}
\end{align*}
are polynomials with respect to $x_i,y_i,z_i,w_i,q_i,p_i$.

Next, let us consider the relation between the polynomial Hamiltonian system \eqref{eq:15} and mKdV hierarchy. In this paper, we can make the birational transformations between the polynomial Hamiltonian system \eqref{eq:15} and mKdV equations.
\begin{theorem}
The birational transformations
\begin{equation}\label{eq:16}
  \left\{
  \begin{aligned}
   x =&q_1,\\
   y =&p_1+2q_1q_2-2p_2p_3+q_3^2+10p_2^3+40q_1p_2q_3+12q_1^2 p_3+60q_1^3q_3+112q_1^2p_2^2+240q_1^4p_2\\
&+120q_1^6-\frac{g}{2},\\
   z =&q_2+2q_1p_3+12q_1p_2^2+10q_1^2 q_3+40q_1^3p_2+24q_1^5,\\
   w =&p_2+q_1^2,\\
   q =&q_3+2q_1p_2+2q_1^3,\\
   p =&p_3+2q_1q_3+8q_1^2p_2+6q_1^4,\\
   s =&2S
   \end{aligned}
  \right. 
\end{equation}
take the Hamiltonian system \eqref{eq:15} to the system
\begin{equation}\label{eq:17}
  \left\{
  \begin{aligned}
   dx =&w dt+(p-6x^2w)dS,\\
   dy =&(70w^2q-gx+42xq^2+56xwp-140x^3 w^2-70x^4 q+20x^7+14x^2z+\alpha_0) dt\\
&+f_1(x,y,z,w,q,p)dS,\\
   dz =&y dt+f_2(x,y,z,w,q,p)dS,\\
   dw =&q dt+f_3(x,y,z,w,q,p)dS,\\
   dq =&p dt+f_4(x,y,z,w,q,p)dS,\\
   dp =&z dt+f_5(x,y,z,w,q,p)dS,
   \end{aligned}
  \right. 
\end{equation}
where $f_i(x,y,z,w,q,p) \in {\mathbb C}[x,y,z,w,q,p] \ (i=1,2,3,4,5)$. 
\end{theorem}
Setting $u:=x$, we see that
\begin{equation}
\frac{\partial u}{\partial t}=w, \quad  \frac{\partial^2 u}{\partial t^2}=q, \quad \frac{\partial^3 u}{\partial t^3}=p, \quad  \frac{\partial^4 u}{\partial t^4}=z, \quad \frac{\partial^5 u}{\partial t^5}=y,
\end{equation}
and
\begin{equation}\label{eq:18}
  \left\{
  \begin{aligned}
   \frac{\partial^6 u}{\partial t^6} =&14u^2\frac{\partial^4 u}{\partial t^4}+56u\frac{\partial u}{\partial t}\frac{\partial^3 u}{\partial t^3}+42u\left(\frac{\partial^2 u}{\partial t^2} \right)^2-70\left(u^4-\left(\frac{\partial u}{\partial t} \right)^2 \right)\frac{\partial^2 u}{\partial t^2}-140u^3\left(\frac{\partial u}{\partial t} \right)^2\\
&+20u^7-gu+\alpha_0,\\
   \frac{\partial u}{\partial S} =&\frac{\partial^3 u}{\partial t^3}-6u^2\frac{\partial u}{\partial t}.
   \end{aligned}
  \right. 
\end{equation}
The first equation in \eqref{eq:18} coincides with an autonomous version of $P_{II}^{(3)}$, and the second equation just coincides with the mKdV equation.

\section{Autonomous version of $P_{II}^{(3)}$ and mKdV5 equation}

In this section, we find a one-parameter family of polynomial Hamiltonian system in two variables given by
\begin{equation}\label{eq:19}
  \left\{
  \begin{aligned}
   dq_1 =&\frac{\partial K_1}{\partial p_1}dt+\frac{\partial K_3}{\partial p_1}ds, \quad dp_1 =-\frac{\partial K_1}{\partial q_1}dt-\frac{\partial K_3}{\partial q_1}ds,\\
   dq_2 =&\frac{\partial K_1}{\partial p_2}dt+\frac{\partial K_3}{\partial p_2}ds, \quad dp_2 =-\frac{\partial K_1}{\partial q_2}dt-\frac{\partial K_3}{\partial q_2}ds,\\
   dq_3 =&\frac{\partial K_1}{\partial p_3}dt+\frac{\partial K_3}{\partial p_3}ds, \quad dp_3 =-\frac{\partial K_1}{\partial q_3}dt-\frac{\partial K_3}{\partial q_3}ds
   \end{aligned}
  \right. 
\end{equation}
with the polynomial Hamiltonians
\begin{align}
\begin{split}
K_1=&q_1^2p_1-\alpha_0 q_1-\frac{p_2^4}{2}-\frac{g}{2}p_2+\frac{1}{2}p_3^2+p_1 p_2-q_2q_3+p_2 q_3^2-2p_2^2 p_3,\\
K_3=&\frac{1}{4}p_1^2-\frac{1}{4}gp_1-\frac{\alpha_0}{2}q_2-\alpha_0 q_1p_3+q_1p_1q_2-\alpha_0 q_1p_2^2+\frac{1}{2}p_1q_3^2-p_1p_2p_3\\
&+q_1^2p_1p_3+q_1^2p_1p_2^2 \quad (g \in {\mathbb C}).
\end{split}
\end{align}
Here, $\alpha_0$ and $\alpha_1$ are constant complex parameters. We remark that the polynomial Hamiltonians $K_2$ and $K_3$ are different (see figure 1). For example, $deg(K_2)=4$ with respect to $q_1,p_1,q_2,p_2,q_3,p_3$. On the other hand, $deg(K_3)=5$ with respect to $q_1,p_1,q_2,p_2,q_3,p_3$.

\begin{figure}
%%%%%%%%%%%\input{PII31}%%%%%%%%%%%%
%WinTpicVersion3.08
\unitlength 0.1in
\begin{picture}( 40.3000, 17.0000)( 11.8000,-28.9000)
% STR 2 0 3 0
% 3 1180 1260 1180 1360 2 0
% Polynomial Hamiltonian
\put(11.8000,-13.6000){\makebox(0,0)[lb]{Polynomial Hamiltonian}}%
% STR 2 0 3 0
% 3 1190 1460 1190 1560 2 0
% system \eqref{eq:15} in two variables
\put(11.9000,-15.6000){\makebox(0,0)[lb]{system \eqref{eq:15} in two variables}}%
% STR 2 0 3 0
% 3 4200 1260 4200 1360 2 0
% Polynomial Hamiltonian
\put(42.0000,-13.6000){\makebox(0,0)[lb]{Polynomial Hamiltonian}}%
% STR 2 0 3 0
% 3 4210 1460 4210 1560 2 0
% system \eqref{eq:19} in two variables
\put(42.1000,-15.6000){\makebox(0,0)[lb]{system \eqref{eq:19} in two variables}}%
% VECTOR 0 0 3 0
% 2 1860 1780 2740 2660
% 
\special{pn 20}%
\special{pa 1860 1780}%
\special{pa 2740 2660}%
\special{fp}%
\special{sh 1}%
\special{pa 2740 2660}%
\special{pa 2708 2600}%
\special{pa 2702 2622}%
\special{pa 2680 2628}%
\special{pa 2740 2660}%
\special{fp}%
% VECTOR 0 0 3 0
% 2 5210 1800 4360 2660
% 
\special{pn 20}%
\special{pa 5210 1800}%
\special{pa 4360 2660}%
\special{fp}%
\special{sh 1}%
\special{pa 4360 2660}%
\special{pa 4422 2628}%
\special{pa 4398 2622}%
\special{pa 4394 2600}%
\special{pa 4360 2660}%
\special{fp}%
% STR 2 0 3 0
% 3 2710 2780 2710 2880 2 0
% Autonomous version of ${P_{II}}^{(3)}$
\put(27.1000,-28.8000){\makebox(0,0)[lb]{Autonomous version of $P_{II}^{(3)}$}}%
% STR 2 0 3 0
% 3 1880 2250 1880 2350 2 0
% $s=0$
\put(18.8000,-23.5000){\makebox(0,0)[lb]{$s=0$}}%
% STR 2 0 3 0
% 3 4840 2250 4840 2350 2 0
% $s=0$
\put(48.4000,-23.5000){\makebox(0,0)[lb]{$s=0$}}%
% STR 2 0 3 0
% 3 1190 1660 1190 1760 2 0
% $(t,s)$
\put(11.9000,-17.6000){\makebox(0,0)[lb]{$(t,s)$}}%
% STR 2 0 3 0
% 3 4200 1650 4200 1750 2 0
% $(t,s)$
\put(42.0000,-17.5000){\makebox(0,0)[lb]{$(t,s)$}}%
% STR 2 0 3 0
% 3 2700 2960 2700 3060 2 0
% in the variable $t$
\put(27.0000,-30.6000){\makebox(0,0)[lb]{in the variable $t$}}%
\end{picture}%
\label{fig:PII31}
\caption{}
\end{figure}

\begin{proposition}
The system \eqref{eq:19} satisfies the compatibility conditions$:$
\begin{align}
\begin{split}
&\frac{\partial }{\partial s} \frac{\partial q_1}{\partial t}=\frac{\partial }{\partial t} \frac{\partial q_1}{\partial s}, \quad \frac{\partial }{\partial s} \frac{\partial p_1}{\partial t}=\frac{\partial }{\partial t} \frac{\partial p_1}{\partial s}, \quad \frac{\partial }{\partial s} \frac{\partial q_2}{\partial t}=\frac{\partial }{\partial t} \frac{\partial q_2}{\partial s},\\
&\frac{\partial }{\partial s} \frac{\partial p_2}{\partial t}=\frac{\partial }{\partial t} \frac{\partial p_2}{\partial s}, \quad \frac{\partial }{\partial s} \frac{\partial q_3}{\partial t}=\frac{\partial }{\partial t} \frac{\partial q_3}{\partial s}, \quad \frac{\partial }{\partial s} \frac{\partial p_3}{\partial t}=\frac{\partial }{\partial t} \frac{\partial p_3}{\partial s}.
\end{split}
\end{align}

\end{proposition}

\begin{proposition}
Two Hamiltonians $K_1$ and $K_3$ satisfy
\begin{equation}
\{K_1,K_3\}=0.
\end{equation}
\end{proposition}

\begin{theorem}
The system \eqref{eq:19} has $K_1$ and $K_3$ as its first integrals.
\end{theorem}

\begin{proposition}
The system \eqref{eq:19} admits the extended affine Weyl group symmetry of type $A_1^{(1)}$ as the group of its B{\"a}cklund transformations, whose generators $s_0,s_1,{\pi}$ defined as follows$:$ with {\it the notation} $(*):=(q_1,p_1,q_2,p_2,q_3,p_3,t,s;\alpha_0,\alpha_1)$\rm{: \rm}
\begin{align*}
s_0:(*) \rightarrow &\left(q_1-\frac{\alpha_0}{f_0},p_1,q_2,p_2,q_3,p_3,t,s;-\alpha_0,\alpha_1+2\alpha_0 \right),\\
s_1:(*) \rightarrow &(q_1-\frac{\alpha_1}{f_1},p_1+\frac{4\alpha_1(q_2+2q_1p_3+2q_1 p_2^2)}{f_1}+\frac{4\alpha_1^2(p_3-p_2^2+4q_1q_3+4q_1^2p_2)}{f_1^2},\\
&q_2-\frac{4\alpha_1(p_3-2q_1^2 p_2)}{f_1}+\frac{8\alpha_1^2 (q_3+q_1p_2-2q_1^3)}{f_1^2}+\frac{8\alpha_1^3(p_2+2q_1^2)}{f_1^3},p_2+\frac{4\alpha_1 q_1}{f_1}-\frac{2\alpha_1^2}{f_1^2},\\
&q_3+\frac{4\alpha_1 (p_2-q_1^2)}{f_1}+\frac{12\alpha_1^2 q_1}{f_1^2}-\frac{4\alpha_1^3}{f_1^3},p_3+\frac{4\alpha_1 q_3}{f_1}+\frac{8\alpha_1^2 (p_2-q_1^2)}{f_1^2}+\frac{16\alpha_1^3 q_1}{f_1^3}-\frac{4\alpha_1^4}{f_1^4},\\
&t,s;\alpha_0+2\alpha_1,-\alpha_1),\\
\pi:(*) \rightarrow &(-q_1,-f_1,-(q_2+4q_1p_3+8q_1^2 q_3+8q_1^3 p_2),-(p_2+2q_1^2),-(q_3+4q_1p_2+4q_1^3),\\
&-(p_3+4q_1q_3+8q_1^2p_2+4q_1^4),t,s;\alpha_1,\alpha_0),
\end{align*}
where $f_0:=p_1$ and $f_1:=p_1+4q_1q_2+2q_3^2-4p_2p_3+4q_1^2 p_3+4q_1^2 p_2^2-g$.
\end{proposition}
Here, the parameters $\alpha_i$ satisfy the relation $\alpha_0+\alpha_1=0$.

\begin{proposition}
The polynomial Hamiltonian system \eqref{eq:19} becomes again a polynomial Hamiltonian system in each coordinate $R_i \ (i=0,1)${\rm : \rm}
\begin{align*}
\begin{split}
R_0:(x_0,y_0,z_0,w_0)=&\left(\frac{1}{q_1},-(q_1f_0-\alpha_0)q_1,q_2,p_2,q_3,p_3 \right),\\
R_1:(x_1,y_1,z_1,w_1)=&(\frac{1}{q_1},-(q_1f_1-\alpha_1)q_1,q_2+4q_1p_3+8q_1^2 q_3+8q_1^3 p_2,p_2+2q_1^2,\\
&q_3+4q_1p_2+4q_1^3,p_3+4q_1q_3+8q_1^2 p_2+4q_1^4 ),
\end{split}
\end{align*}
where $f_0:=p_1$ and $f_1:=p_1+4q_1q_2+2q_3^2-4p_2p_3+4q_1^2 p_3+4q_1^2 p_2^2-g$.
\end{proposition}

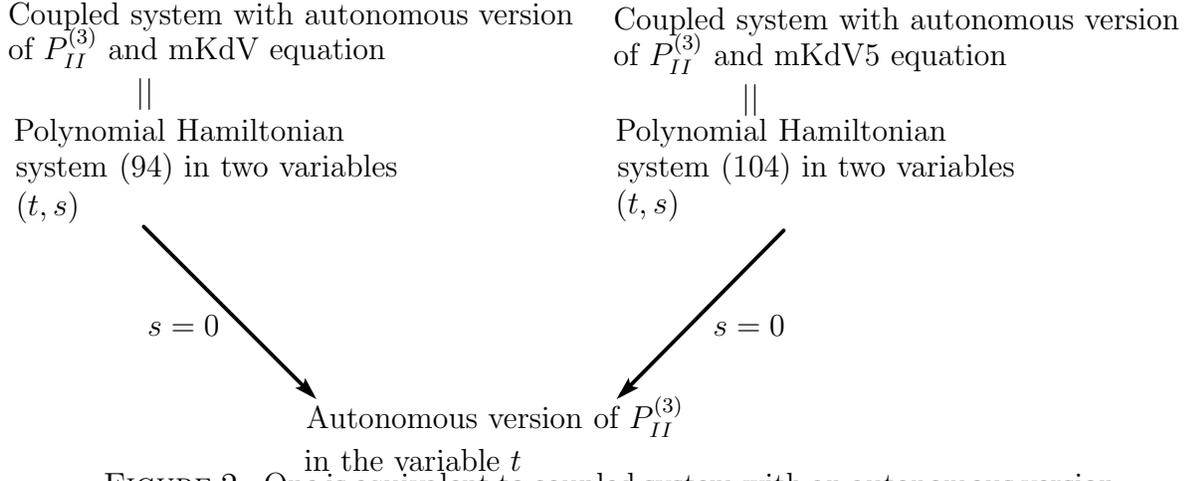
\begin{figure}
\begin{flushleft}
%%%%%%%%%%%%\input{PII32}%%%%%%%%%%%
%WinTpicVersion3.08
\unitlength 0.1in
\begin{picture}( 40.6000, 23.1000)( 11.5000,-28.9000)
% STR 2 0 3 0
% 3 1180 1260 1180 1360 2 0
% Polynomial Hamiltonian
\put(11.8000,-13.6000){\makebox(0,0)[lb]{Polynomial Hamiltonian}}%
% STR 2 0 3 0
% 3 1190 1460 1190 1560 2 0
% system \eqref{eq:15} in two variables
\put(11.9000,-15.6000){\makebox(0,0)[lb]{system \eqref{eq:15} in two variables}}%
% STR 2 0 3 0
% 3 4330 1260 4330 1360 2 0
% Polynomial Hamiltonian
\put(43.3000,-13.6000){\makebox(0,0)[lb]{Polynomial Hamiltonian}}%
% STR 2 0 3 0
% 3 4340 1460 4340 1560 2 0
% system \eqref{eq:19} in two variables
\put(43.4000,-15.6000){\makebox(0,0)[lb]{system \eqref{eq:19} in two variables}}%
% VECTOR 0 0 3 0
% 2 1860 1780 2740 2660
% 
\special{pn 20}%
\special{pa 1860 1780}%
\special{pa 2740 2660}%
\special{fp}%
\special{sh 1}%
\special{pa 2740 2660}%
\special{pa 2708 2600}%
\special{pa 2702 2622}%
\special{pa 2680 2628}%
\special{pa 2740 2660}%
\special{fp}%
% VECTOR 0 0 3 0
% 2 5210 1800 4360 2660
% 
\special{pn 20}%
\special{pa 5210 1800}%
\special{pa 4360 2660}%
\special{fp}%
\special{sh 1}%
\special{pa 4360 2660}%
\special{pa 4422 2628}%
\special{pa 4398 2622}%
\special{pa 4394 2600}%
\special{pa 4360 2660}%
\special{fp}%
% STR 2 0 3 0
% 3 2710 2780 2710 2880 2 0
% Autonomous version of ${P_{II}}^{(3)}$
\put(27.1000,-28.8000){\makebox(0,0)[lb]{Autonomous version of $P_{II}^{(3)}$}}%
% STR 2 0 3 0
% 3 1880 2250 1880 2350 2 0
% $s=0$
\put(18.8000,-23.5000){\makebox(0,0)[lb]{$s=0$}}%
% STR 2 0 3 0
% 3 4840 2250 4840 2350 2 0
% $s=0$
\put(48.4000,-23.5000){\makebox(0,0)[lb]{$s=0$}}%
% STR 2 0 3 0
% 3 1190 1660 1190 1760 2 0
% $(t,s)$
\put(11.9000,-17.6000){\makebox(0,0)[lb]{$(t,s)$}}%
% STR 2 0 3 0
% 3 4330 1650 4330 1750 2 0
% $(t,s)$
\put(43.3000,-17.5000){\makebox(0,0)[lb]{$(t,s)$}}%
% STR 2 0 3 0
% 3 2700 2960 2700 3060 2 0
% in the variable $t$
\put(27.0000,-30.6000){\makebox(0,0)[lb]{in the variable $t$}}%
% STR 2 0 3 0
% 3 1820 1070 1820 1170 2 0
% $||$
\put(18.2000,-11.7000){\makebox(0,0)[lb]{$||$}}%
% STR 2 0 3 0
% 3 1150 650 1150 750 2 0
% Coupled system with autonomous version
\put(11.5000,-7.5000){\makebox(0,0)[lb]{Coupled system with autonomous version}}%
% STR 2 0 3 0
% 3 1150 850 1150 950 2 0
% of ${P_{II}}^{(3)}$ and mKdV equation
\put(11.5000,-9.5000){\makebox(0,0)[lb]{of $P_{II}^{(3)}$ and mKdV equation}}%
% STR 2 0 3 0
% 3 4990 1100 4990 1200 2 0
% $||$
\put(49.9000,-12.0000){\makebox(0,0)[lb]{$||$}}%
% STR 2 0 3 0
% 3 4320 680 4320 780 2 0
% Coupled system with autonomous version
\put(43.2000,-7.8000){\makebox(0,0)[lb]{Coupled system with autonomous version}}%
% STR 2 0 3 0
% 3 4320 880 4320 980 2 0
% of ${P_{II}}^{(3)}$ and mKdV5 equation
\put(43.2000,-9.8000){\makebox(0,0)[lb]{of $P_{II}^{(3)}$ and mKdV5 equation}}%
\end{picture}%
\label{fig:PII32}
\caption{One is equivalent to coupled system with an autonomous version of $P_{II}^{(3)}$ and mKdV equation and the other is equivalent to coupled system with an autonomous version of $P_{II}^{(3)}$ and mKdV5 equation.}
\end{flushleft}
\end{figure}

Next, let us consider the relation between the polynomial Hamiltonian system \eqref{eq:19} and mKdV5 equation.
\begin{theorem}
The birational transformations
\begin{equation}\label{eq:20}
  \left\{
  \begin{aligned}
   x =&q_1,\\
   y =&p_1+2q_1q_2-2p_2p_3+q_3^2+10p_2^3+40q_1p_2q_3+12q_1^2 p_3+60q_1^3q_3+112q_1^2p_2^2+240q_1^4p_2\\
&+120q_1^6-\frac{g}{2},\\
   z =&q_2+2q_1p_3+12q_1p_2^2+10q_1^2 q_3+40q_1^3p_2+24q_1^5,\\
   w =&p_2+q_1^2,\\
   q =&q_3+2q_1p_2+2q_1^3,\\
   p =&p_3+2q_1q_3+8q_1^2p_2+6q_1^4,\\
   s =&2S
   \end{aligned}
  \right. 
\end{equation}
take the Hamiltonian system \eqref{eq:19} to the system
\begin{equation}\label{eq:21}
  \left\{
  \begin{aligned}
   dx =&w dt+(y-10x^2 p-40xwq-10w^3+30x^4 w)dS,\\
   dy =&(70w^2q-gx+42xq^2+56xwp-140x^3 w^2-70x^4 q+20x^7+14x^2z+\alpha_0) dt\\
&+g_1(x,y,z,w,q,p)dS,\\
   dz =&y dt+g_2(x,y,z,w,q,p)dS,\\
   dw =&q dt+g_3(x,y,z,w,q,p)dS,\\
   dq =&p dt+g_4(x,y,z,w,q,p)dS,\\
   dp =&z dt+g_5(x,y,z,w,q,p)dS,
   \end{aligned}
  \right. 
\end{equation}
where $g_i(x,y,z,w,q,p) \in {\mathbb C}[x,y,z,w,q,p] \ (i=1,2,3,4,5)$. 
\end{theorem}
Setting $u:=x$, we see that
\begin{equation}
\frac{\partial u}{\partial t}=w, \quad  \frac{\partial^2 u}{\partial t^2}=q, \quad \frac{\partial^3 u}{\partial t^3}=p, \quad  \frac{\partial^4 u}{\partial t^4}=z, \quad \frac{\partial^5 u}{\partial t^5}=y,
\end{equation}
and
\begin{equation}\label{eq:22}
  \left\{
  \begin{aligned}
   \frac{\partial^6 u}{\partial t^6} =&14u^2\frac{\partial^4 u}{\partial t^4}+56u\frac{\partial u}{\partial t}\frac{\partial^3 u}{\partial t^3}+42u\left(\frac{\partial^2 u}{\partial t^2} \right)^2-70\left(u^4-\left(\frac{\partial u}{\partial t} \right)^2 \right)\frac{\partial^2 u}{\partial t^2}-140u^3\left(\frac{\partial u}{\partial t} \right)^2\\
&+20u^7-gu+\alpha_0,\\
   \frac{\partial u}{\partial S} =&\frac{\partial^5 u}{\partial t^5}-\left( 10u^2 \frac{\partial^3 u}{\partial t^3}+40u\frac{\partial u}{\partial t} \frac{\partial^2 u}{\partial t^2}+10\left(\frac{\partial u}{\partial t} \right)^3-30u^4 \frac{\partial u}{\partial t} \right).
   \end{aligned}
  \right. 
\end{equation}
The first equation in \eqref{eq:22} coincides with an autonomous version of $P_{II}^{(3)}$, and the second equation just coincides with mKdV5 equation.

\section{Autonomous version of the second Painlev\'e system and the mKdV equation}

In this section, we study a one-parameter family of partial differential systems of polynomial type in two variables given by
\begin{equation}\label{eq:A11}
  \left\{
  \begin{aligned}
   dq_1 =&\frac{\partial q_1}{\partial t}dt+\frac{\partial q_1}{\partial s}ds\\
   =&\left(q_1^2+p_1-\frac{1}{2} \right)dt+\left(q_1^2+p_1-\frac{1}{2} \right)ds,\\
   dp_1 =&\frac{\partial p_1}{\partial t}dt+\frac{\partial p_1}{\partial s}ds\\
   =&(-2q_1p_1-\alpha)dt+(-2q_1p_1-\alpha)ds.
   \end{aligned}
  \right. 
\end{equation}
We easily see that the system \eqref{eq:A11} satisfies the compatibility conditions$:$
\begin{equation}
\frac{\partial }{\partial s} \frac{\partial q_1}{\partial t}=\frac{\partial }{\partial t} \frac{\partial q_1}{\partial s}, \quad \frac{\partial }{\partial s} \frac{\partial p_1}{\partial t}=\frac{\partial }{\partial t} \frac{\partial p_1}{\partial s}.
\end{equation}
In this paper, we will make the birational transformations between the system \eqref{eq:A11} and the mKdV equation
\begin{equation}
   \frac{\partial^3 u}{\partial t^3} =6u^2 \frac{\partial u}{\partial t}-\frac{\partial u}{\partial s}.
\end{equation}

\begin{proposition}
The system \eqref{eq:A11} has the Hamiltonian $H$
\begin{equation}
H=q_1^2 p_1+\frac{p_1^2}{2}-\frac{p_1}{2}+\alpha q_1
\end{equation}
as its first integrals.
\end{proposition}

\begin{proposition}
The system \eqref{eq:A11} admits the extended affine Weyl group symmetry of type $A_1^{(1)}$ as the group of its B{\"a}cklund transformations, whose generators $s_0,s_1,{\pi}$ defined as follows$:$ with {\it the notation} $(*):=(q_1,p_1,t,s;\alpha)$\rm{; \rm}
\begin{align*}
s_0:(*) \rightarrow &\left(q_1+\frac{\alpha}{f_0},p_1,t,s;-\alpha \right),\\
s_1:(*) \rightarrow &\left(q_1-\frac{\alpha}{f_1},p_1+\frac{4\alpha q_1}{f_1}-\frac{2\alpha^2}{f_1^2},t,s;-\alpha \right),\\
\pi:(*) \rightarrow &(-q_1,-f_1,t,s;-\alpha),
\end{align*}
where $f_0:=p_1$ and $f_1:=p_1+2q_1^2-1$.
\end{proposition}

\begin{proposition}
The system \eqref{eq:A11} admits some particular solutions$:$
\begin{equation}
(q_1,p_1;\alpha)=\left(0,\frac{1}{2};0 \right)
\end{equation}
and
\begin{equation}
(q_1,p_1;\alpha)=\left(\pm \frac{1}{\sqrt{2}},0;0 \right),
\end{equation}
and
\begin{equation}
(q_1,p_1;\alpha)=\left(-\frac{1}{\sqrt{2}} \rm{tanh \rm} \left(\frac{t+s+c}{\sqrt{2}} \right),0;0 \right) \quad (c \in {\mathbb C}).
\end{equation}
\end{proposition}

\begin{theorem}
The system \eqref{eq:A11} is invariant under the following auto-B{\"a}cklund transformations $T_0,T_1$ defined as follows$:$ with {\it the notation} $(*):=(q_1,p_1,t,s;\alpha)$\rm{; \rm}
\begin{align*}
T_0:(*) \rightarrow &\left(-q_1+\frac{\alpha}{p_1+2q_1^2-1},1-2q_1^2-p_1,t,s;\alpha \right),\\
T_1:(*) \rightarrow &\left(-q_1-\frac{\alpha}{p_1},1-p_1-\frac{2(q_1p_1+\alpha)^2}{p_1^2},t,s;\alpha \right),
\end{align*}
where $T_0:=\pi s_0, \quad T_1:=\pi s_1$.
\end{theorem}
Applying these B{\"a}cklund transformations $T_0^m$ and $T_1^n$ $(m,n=1,2,3,..)$, we can obtain a series of its particular solutions.

Next, let us consider the relation between the system \eqref{eq:A11} and the mKdV equation. In this paper, we can make the birational transformations between the system \eqref{eq:A11} and the mKdV equation.
\begin{theorem}
The birational transformations
\begin{equation}\label{eq:A12}
  \left\{
  \begin{aligned}
   x =&q_1,\\
   y =&q_1^2+p_1-\frac{1}{2}
   \end{aligned}
  \right. 
\end{equation}
take the system \eqref{eq:A11} to the system
\begin{equation}\label{eq:A13}
  \left\{
  \begin{aligned}
   dx =&y dt+y ds,\\
   dy =&(2x^3-x-\alpha) dt+(2x^3-x-\alpha)ds.
   \end{aligned}
  \right. 
\end{equation}
\end{theorem}

Setting $u:=x$, we see that
\begin{equation}
\frac{\partial u}{\partial t}=\frac{\partial u}{\partial s}=y
\end{equation}
and
\begin{equation}\label{eq:A14}
  \left\{
  \begin{aligned}
   \frac{\partial^2 u}{\partial t^2} =&2u^3-u-\alpha,\\
   \frac{\partial u}{\partial s} =&\frac{\partial u}{\partial t}.
   \end{aligned}
  \right. 
\end{equation}
The first equation in \eqref{eq:A14} coincides with an autonomous version of $P_{II}$.

Making a partial derivation in the variable $t$ for the first equation of \eqref{eq:A14}, we obtain
\begin{equation}\label{eq:A15}
  \left\{
  \begin{aligned}
   \frac{\partial^3 u}{\partial t^3} =&6u^2 \frac{\partial u}{\partial t}-\frac{\partial u}{\partial t},\\
   0=&\frac{\partial u}{\partial t}-\frac{\partial u}{\partial s}.
   \end{aligned}
  \right. 
\end{equation}

Adding each system in \eqref{eq:A15}, we can obtain the mKdV equation
\begin{equation}\label{eq:A16}
   \frac{\partial^3 u}{\partial t^3} =6u^2 \frac{\partial u}{\partial t}-\frac{\partial u}{\partial s}.
\end{equation}
It is known that the Hamiltonian system with a special parameter $\alpha=0$
\begin{equation}\label{eq:ABB}
  \left\{
  \begin{aligned}
   \frac{dq_1}{dt} =&\frac{\partial H}{\partial p_1}=q_1^2+p_1-\frac{1}{2},\\
   \frac{dp_1}{dt} =&-\frac{\partial H}{\partial q_1}=-2q_1p_1
   \end{aligned}
  \right. 
\end{equation}
can be solved by using Jacobi' s elliptic function $sn(t)$. Setting its solution $\varphi(t)$
$$
\frac{d^2 \varphi(t)}{d t^2}=2\varphi(t)^3-\varphi(t),
$$
we can make a stationary solution $u(t, s)$ of the system \eqref{eq:A16};
\begin{equation}
u(t, s) := \varphi(t + s + c) \quad (c : constant).
\end{equation}

{\it Acknowledgements.} The author would like to thank Prof. S. Iwao, K. Fuji, Y Ohta, Y. Ohyama, N. Suzuki, H. Watanabe for useful discussions, and Prof. M. Jimbo  and S. Kakei gave helpful advice and encouragement. This work is supported by Kyoto university, Aoyama Gakuin  university and Kitami Institute of Technology.


\begin{thebibliography}{99}
\bibitem[1]{Clarkson} P. Clarkson, N. Joshi and A. Pickering, 
{\em B{\"a}cklund transformations for the second Painlev\'e hierarchy: a modified truncation approach}, 
Inverse Problem. {\bf 15} (1999), 175--187.



\bibitem[2]{Mazzocco} M. Mazzocco and M. Y. Mo, 
{\em The Hamiltonian structure of the second Painlev\'e hierarchy}, 
http://arXiv.org/abs/nlin.SI/0610066. 

 

\bibitem[3]{Joshi} N. Joshi, 
{\em The second Painlev\'e hierarchy and the stationary KdV hierarchy}, 
Publ. RIMS. Kyoto Univ. {\bf 40} (2004), 1039-1061. 


\bibitem[4]{Sasano1} Y. Sasano,
{\em Symmetry and holomorphy of the second member of the second Painlev\'e hierarchy}, preprint.

\bibitem[5]{5} P. D. Lax, 
{\em Integrals of nonlinear equations of evolution and solitary waves}, Commun. Pure Appl. Math. {\bf 21} (1968), 467--490. 


\bibitem[6]{6} P. Painlev\'e, {\em M\'emoire sur les \'equations diff\'erentielles dont l'int\'egrale g\'en\'erale est uniforme}, Bull. Soci\'et\'e Math\'ematique de France. {\bf 28} (1900),  201--261.

\bibitem[7]{7} P. Painlev\'e, {\em Sur les \'equations diff\'erentielles du second ordre et d'ordre sup\'erieur dont l'int\'egrale est uniforme}, Acta Math. {\bf 25} (1902), 1--85.

\bibitem[8]{Cos1} C. M. Cosgrove, 
{\em Chazy classes IX-XI of third-order differential equations}, 
Stud. Appl. Math. {\bf 104}, (2000), 171--228. 

\bibitem[9]{9} C. M. Cosgrove, 
{\em Higher-order Painlev\'e equations in the polynomial class I. Bureau symbol P2}, Stud. Appl. Math. {\bf 104}, (2000), 1--65.

\bibitem[10]{Bureau} F. Bureau, 
{\em Differential equations with fixed critical points}, 
Anna Matemat. {\bf 66} (1964), 1--116.

\bibitem[11]{11} J. Chazy, 
{\em Sur les \'equations diff\'erentielles dont l'int\'egrale g\'en\'erale est uniforme et admet des singularit\'es essentielles mobiles}, 
Comptes Rendus de l'Acad\'emie des Sciences, Paris. {\bf 149} (1909), 563--565. 

\bibitem[12]{12} J. Chazy, 
{\em Sur les \'equations diff\'erentielles dont l'int\'egrale g\'en\'erale poss\'ede une coupure essentielle mobile }, 
Comptes Rendus de l'Acad\'emie des Sciences, Paris. {\bf 150} (1910), 456--458. 


\bibitem[13]{13} J. Chazy, 
{\em Sur les \'equations diff\'erentielles du trousi\'eme ordre et d'ordre sup\'erieur dont l'int\'egrale a ses points critiques fixes}, 
Acta Math. {\bf 34} (1911), 317--385.

\bibitem[14]{Cosgrove1} C. M. Cosgrove and G. Scoufis,
{\em Painlev\'e classification of a class of differential equations of the second order and second degree}, Studies in Applied Mathematics. {\bf 88} (1993), 25-87.

\bibitem[15]{Cosgrove2} C. M. Cosgrove,
{\em All binomial-type Painlev\'e equations of the second order and degree three or higher}, Studies in Applied Mathematics. {\bf 90} (1993), 119-187.

\bibitem[16]{Ince} E. L. Ince, 
{\em Ordinary differential equations}, Dover Publications, New York, (1956).

\bibitem[17]{FS1} K. Fuji and T. Suzuki, {\em The sixth Painlev\'e equation arising from $D_4^{(1)}$ hierarchy}, J. Phys. A: Math. Gen. {\bf 39} (2006), 12073--12082.

\bibitem[18]{FS2} K. Fuji and T Suzuki, {\em Drinfeld-Sokolov hierarchies of type A and fourth order Painlev\'e systems}, Funkcial. Ekvac {\bf 53} (2010), 143--167.


\bibitem[20]{N1} M. Noumi and Y. Yamada, 
        {\em Higher order Painlev\'e equations of type $A_l^{(1)}$}, Funkcial. Ekvac. {\bf 41} (1998), 483--503.

\bibitem[21]{N2} M. Noumi and Y. Yamada,
        {\em Affine Weyl Groups, Discrete Dynamical Systems and Painlev\`e Equations}, Comm Math Phys {\bf 199} (1998), 281--295.

\bibitem[22]{N3} M. Noumi and Y. Yamada,
        {\em A new Lax pair for the sixth Painlev\'e equation associated with ${\mathfrak so}(8)$}, in Microlocal Analysis and Complex Fourier Analysis. ed. T.Kawai and K.Fujita, (World Scientific, 2002), 238--252.


\bibitem[23]{N4} M. Noumi,
        {\em Painlev\'e Equations Through Symmetry}, American Mathematical Society, {\bf 223} (2004).





\end{thebibliography}
\end{document}